 \documentclass[final,authoryear,3p,times,fleqn]{elsarticle}

\usepackage{amssymb}
\usepackage{amsmath}
\usepackage{amsthm}
\usepackage[ruled]{algorithm2e}

\usepackage{graphicx}
\usepackage{float} 
\usepackage{subcaption}
\usepackage{pxfonts,times}
\usepackage{siunitx}

\usepackage{bm}
\usepackage{color}
\usepackage{comment}
\usepackage{natbib}

\usepackage[
	colorlinks=true, 
	citecolor=blue, 
	urlcolor=blue]{hyperref}

\newtheorem{theo}{Theorem}
\newtheorem{defi}{Definition}
\newtheorem{lemm}{Lemma}
\newtheorem{prop}{Proposition}

\newtheorem{coro}{Corollary}

\newtheorem*{prf}{Proof}

\makeatletter
\def\th@plain{\upshape}
\makeatother

\makeatletter
	
	\@addtoreset{equation}{section}
\makeatother

\journal{Elsevier}

\begin{document}

\begin{frontmatter}



\title{Stability analysis of a departure time choice problem with atomic vehicle models}



\author[a]{Koki Satsukawa\corref{cor1}}
\address[a]{Institute of Transdisciplinary Sciences for Innovation, Kanazawa University, Ishikawa, Japan}
\ead{satsukawa@staff.kanazawa-u.ac.jp} 
\cortext[cor1]{Corresponding author.}

\author[b]{Kentaro Wada}
\address[b]{Institute of Systems and Information Engineering, University of Tsukuba, Ibaraki, Japan}
\ead{wadaken@sk.tsukuba.ac.jp}

\author[c]{Takamasa Iryo}
\address[c]{Graduate School of Information Sciences, Tohoku University, Miyagi, Japan}
\ead{iryo@tohoku.ac.jp}

\begin{abstract}
\color{black}
In this study, we analyse the global stability of the equilibrium in a departure time choice problem using a game-theoretic approach that deals with atomic users. 
We first formulate the departure time choice problem as a strategic game in which atomic users select departure times to minimise their trip cost; we call this game the `departure time choice game'.
The concept of the epsilon-Nash equilibrium is introduced to ensure the existence of pure-strategy equilibrium corresponding to the departure time choice equilibrium in conventional fluid models.
Then, we prove that the departure time choice game is a weakly acyclic game.
By analysing the convergent better responses, we clarify the mechanisms of global convergence to equilibrium. 
This means that the epsilon-Nash equilibrium is achieved by sequential better responses of users, which are departure time changes to improve their own utility, in an appropriate order.
Specifically, the following behavioural rules are important to ensure global convergence: (i) the adjustment of the departure time of the first user departing from the origin to the corresponding equilibrium departure time and (ii) the fixation of users to their equilibrium departure times in order (starting with the earliest). 
Using convergence mechanisms, we construct evolutionary dynamics under which global stability is guaranteed. 
We also investigate the stable and unstable dynamics studied in the literature based on convergence mechanisms, and gain insight into the factors influencing the different stability results.
Finally, numerical experiments are conducted to demonstrate the theoretical results.
\color{black}
\end{abstract}


\end{frontmatter}



\section{Introduction}

\color{black}
\subsection{Background and purpose}\label{Sec:Backgound}
Stability is an important property that ensures the realisation of the equilibrium in transport systems.
It has been discussed for decades since the seminal work of \cite{Beckmann1956-vr}, who indicated the significance of analysing stability by stating that `\textit{An equilibrium would be just an extreme state of rare occurrence if it were not stable---that is, if there were no forces which tended to restore equilibrium as soon as small deviations from it occurred}'.
Such a force leading a traffic state towards equilibrium has been studied by analysing the day-to-day adjustment behaviour of users' choices, which is called \textit{evolutionary dynamics or day-to-day dynamics}.
\cite{Smith1984-ed} developed Smith dynamics~\citep[named in][]{Hofbauer2009-wr}, and proved the asymptotic stability of user equilibrium in a conventional static traffic assignment problem.
\cite{Friesz1994-mb,Zhang1996-ir} proved that stability holds under different evolutionary dynamics in static traffic assignment problems.

However, stability is not necessarily a general property in dynamic traffic assignment problems~\citep[][]{Iryo2013-pe}.
One of the representative problems is the departure time choice (DTC) problem in a single-bottleneck network proposed by \cite{Vickrey1969-rg}.
The DTC problem represents the within-day traffic flow dynamics (formation and dissipation of queues) in a simple and tractable manner using a bottleneck (i.e. point queue) model and describes the equilibrium distribution of users' departure times.
In this problem, both positive and negative stability results have been presented.
\cite{Iryo2008-ep,Guo2018-fy,Iryo2019-co} analysed the behaviour of rational evolutionary dynamics in the sense that commuters prefer departure times with lower trip costs, i.e. they conduct better responses.
These studies demonstrated the instability of the equilibrium and its non-convergence under the evolutionary dynamics.
\cite{Jin2021-xe} proposed \textit{local} evolutionary dynamics, under which users can change their departure times only to their adjacent ones within the equilibrium rush-hour period.
The author then proved the local asymptotic stability of the equilibrium within a certain stability region.

Considering these results, there remains room for further studies on the stability analysis of DTC problems in the following areas.
First, little is known about the theoretical results of the global stability of the equilibrium.
In particular, it is important to investigate whether equilibrium can be reached from an arbitrary initial traffic state (i.e. global convergence).\footnote{\cite{Jin2020-tz} showed the globally stable evolutionary dynamics. However, this study assumes that users select both origin departure and destination arrival times, which is different from the DTC problem focused on in this study.}
Furthermore, it is worthwhile to clarify the mechanisms underlying the different stability results demonstrated in the literature.
Previous studies first selected the evolutionary dynamics to be analysed and investigated the corresponding behaviour of the traffic flows, which helped to clarify the properties of each dynamics separately.
Meanwhile, it is not well understood which differences in the behavioural rules between stable and unstable dynamics are the main factors leading to the different stability results in the DTC problems.

To tackle these research gaps, we analyse the stability of equilibrium in DTC problems using \textit{an atomic approach} that deals with atomic users. 
Although the concept of continuous traffic flow in a conventional fluid approach is common and useful for theoretical analysis, there are merits to this atomic approach.
First, it is natural to represent traffic flow by atomic users since it consists of unsplittable individuals in the real world.
Indeed, several studies have analysed traffic assignment problems involving atomic users~\citep[e.g.][]{Rosenthal1973-gz,Iryo2010-zw,Scarsini2018-kb,Cao2021-pg}.
Furthermore, the atomic representation allows us to describe the adjustment behaviour of each individual user.
In the fluid approach, we need to aggregate the adjustment behaviour to formulate the evolutionary dynamics of continuous traffic flow, which approximates a microscopic specification of the individual behaviour~\citep[][]{Sandholm2003-bv,Iryo2016-pj}.
As the atomic approach does not require such aggregation, it can explicitly describe how each individual user changes his/her strategy under different behavioural rules and how the traffic state changes according to these strategy changes.

Moreover, the atomic approach has a significant advantage when incorporating the \textit{ordering property}, which is useful for analysing dynamic traffic assignment problems.
The ordering property is the existence of the assignment order such that equilibrium is achieved by changing the strategy of each user to the equilibrium strategy one by one in the order.
This property was originally proposed by \cite{Kuwahara1993-xm} for the theoretical analysis of the dynamic route-choice equilibrium in a single-origin network;
they showed that the ordering property holds in the network and proved the existence of the equilibrium.
Several studies have used it in combination with fluid approaches for the equilibrium analysis of the route-choice problem in such a special class of networks~\citep[][]{Akamatsu2003-jy,Iryo2018-it,Wada2019-tu}.

The ordering property has been shown to be a powerful tool in disequilibrium analysis (e.g. convergence and stability) of dynamic traffic assignment with route choices when combined with the atomic approach.
For example, \cite{Iryo2011-bk} showed that it is easy to sort and assign individual atomic users according to the assignment order when the ordering property holds, and established an exact and convergent solution algorithm in a single-origin or single-destination network.
Furthermore, \cite{Satsukawa2019-lq} extended the exact solution algorithm to a new method for disequilibrium analysis and succeeded in proving the convergence and stability of the equilibrium.
Apart from these studies, \cite{Waller2006-vx} implicitly assumed the atomic representation of traffic flow and showed the basic framework of an algorithm to compute the exact equilibrium in a single-destination network, although the convergence was not formally proved.\footnote{\textcolor{black}{\cite{Kuwahara1993-xm,Akamatsu2000-zj} presented the solution algorithms for the dynamic route-choice equilibrium using fluid approaches. However, they did not succeed in proving the convergence because it is not easy to sort splittable continuous traffic flow according to the assignment order, which makes it difficult to obtain rigorous results on convergence and stability.}}
Therefore, we use the atomic approach for analysing the convergence and stability of the equilibrium in DTC problems, as a different approach from the conventional fluid one.

In this study, we first formulate a DTC problem as a strategic game in which atomic users select departure times to minimise their trip costs (we call this the `DTC game').
A DTC game is a class of games in which existing studies have pointed out the non-existence of pure-strategy Nash equilibrium~\citep[e.g.][]{Otsubo2008-qj,Ziegelmeyer2008-hi,Ameli2022-zm}.
To overcome this issue, we introduce the concept of the epsilon-Nash equilibrium, which is a weakened form of the pure-strategy Nash equilibrium.
We show its existence and the correspondence between the epsilon-Nash equilibrium and the equilibrium in DTC problems which use a conventional fluid approach.

We then analyse the global convergence to the epsilon-Nash equilibrium.
First, we prove that the DTC game is a weakly acyclic game based on the ordering property.
This means that there exists a better response path (i.e. a sequence of users' better responses) that converges to the equilibrium from any initial traffic state in the DTC game.
Such a convergent better response path provides us with the mechanisms behind the convergence, i.e. the user behaviour necessary to realise the convergent path.
Using these mechanisms, we develop evolutionary dynamics that guarantees convergence to the equilibrium.
Finally, we investigate the reason for the different stability results of the existing evolutionary dynamics based on the convergence mechanisms.

The contributions of this study are summarised as follows.
First, we have extended the atomic approach to the stability analysis of DTC problems by combining the following concepts: the epsilon-Nash equilibrium, weakly acyclic games in game theory and the ordering property of the bottleneck model.
Second, the atomic approach revealed the existence of a better response path that converges globally to the equilibrium based on the ordering property.
Such a convergent path showed the following behavioural rules are important to ensure global convergence: (i) the adjustment of the earliest departing user to the earliest equilibrium departure time, and (ii) the fixation of the departure times of equilibrated users.
Third, we constructed evolutionary dynamics under which global convergence is guaranteed based on the convergence mechanisms.
We also gained insight into the factors influencing the different results of the existing stable and unstable dynamics by examining whether they satisfy the convergence mechanisms.

\color{black}

The remainder of this paper is organised as follows. 
\textcolor{black}{The rest of this section provides a literature review.}
In Section 2, we define the DTC game. 
Section 3 presents an equilibrium analysis of the DTC game while introducing the epsilon-Nash equilibrium.
Section 4 proves that the DTC game is a weakly acyclic game by utilising the ordering property.
We then examine the convergence mechanisms and propose the evolutionary dynamics under which global convergence is guaranteed.
Section 5 presents the numerical experiments to demonstrate the theoretical results.
Finally, Section 6 concludes the paper.

\color{black}
\subsection{Literature review}

This study conducts a disequilibrium analysis of dynamic user equilibrium (DUE) in DTC problems.
Therefore, we provide a literature review of the following research areas related to this study: the equilibrium analysis of DTC problems with a fluid approach, equilibrium analysis with an atomic approach, disequilibrium analysis with a fluid approach and the disequilibrium analysis of other types of DUE problems, such as route-choice equilibrium problems.
For a more comprehensive review of dynamic route-choice and DTC equilibrium, see \cite{Iryo2013-pe} and \cite{Li2020-tk}.

Since the seminal work of \cite{Vickrey1969-rg}, DTC problems have attracted significant attention as the basis for analysing rush-hour traffic congestion
Several studies have conducted equilibrium analyses, such as existence, uniqueness and the closed-form analytical solutions~\citep[][]{Hendrickson1981-bq,Smith1984-ml,Daganzo1985-ko,Arnott1994-ss,Lindsey2004-eu}.
While they considered the heterogeneity of the desired arrival time, several studies also considered the heterogeneity of the value of time~\citep[][]{Arnott1994-ss,Liu2015-xc,Takayama2017-ge}.
Moreover, some studies developed mathematical programming approaches that shows the equivalence to linear programming~\citep[][]{Iryo2005-fo,Iryo2007-fq,Akamatsu2021-ys}.
However, as they focused on the properties of traffic flow and trip costs at equilibrium, there is a lack of understanding of the dynamical process of equilibrium, such as how equilibrium is achieved and whether it is stable.

Some studies have analysed the equilibrium of DTC problems using an atomic approach. 
\cite{Levinson2005-jy} formulated a strategic game in which two or three players choose one of several discrete departure times.
\cite{Zou2006-kl} extended this study to a more general $n$-player strategic game.
This study analysed the pure-strategy Nash equilibrium by numerical experiments and showed that the equilibrium may not exist.
\cite{Ziegelmeyer2008-hi} introduced the concept of mixed-strategy equilibrium as the equilibrium of a strategic game for laboratory studies that assessed the information effect of past departure rates on congestion levels and travel costs.
All the above studies assumed that only one user per unit of time can exit the bottleneck.
For this issue, \cite{Otsubo2008-qj} formulated a DTC game under an arbitrary capacity setting and obtained theoretical results such as the existence and uniqueness of the mixed-strategy equilibrium.
\cite{Ameli2022-zm} have recently proposed a mean-field game approach and formulated a DTC problem using the bathtub model~\citep[][]{Arnott2013-oq} with users whose desired arrival times were different.
Therefore, the application of the atomic approach to more complex situations is currently being explored.
However, no studies have used an atomic approach for the disequilibrium analysis.

Few studies have been conducted on the disequilibrium analysis of DTC problems, and almost all relevant studies have focused on the classical framework with a single bottleneck and homogeneous users.
\cite{Iryo2008-ep} studied the behaviour of Smith dynamics and showed that the dynamics may fail to converge to equilibrium because of the non-monotonicity of the schedule delay cost function.
\cite{Guo2018-fy} analysed several evolutionary dynamics whose adjustment process is similar to Smith dynamics in the sense that commuters prefer departure times with lower trip costs.
This study showed the non-convergence when the number of users exceeded a certain threshold.
In addition, \cite{Iryo2019-co} proved that the equilibrium is not Lyapunov stable under the replicator dynamics.
However, \cite{Jin2021-xe} has recently succeeded in establishing stable evolutionary dynamics for the DTC problem.
This study proposed the local evolutionary dynamics inspired by \cite{Vickrey2020-dt}.
Under the local dynamics, the equilibrium becomes locally asymptotically stable; 
a traffic state converges to the equilibrium if the departure times of traffic flows are within those of the equilibrium.
Although these studies are important building blocks for investigating the stability of the equilibrium, the key factors that determine the stability and instability of the equilibrium are not yet clear.

In addition to the DTC problems, several studies have conducted disequilibrium analyses of DUE problems with route choices.
They have proved stability of the equilibrium in limited classes of network structures.
\cite{Mounce2006-nx} showed the globally asymptotic stability under Smith dynamics in a network in which each route contains only one bottleneck using the Lyapunov approach.
This study proved the existence of the Lyapunov function based on the monotonicity of the route travel time function, which was originally shown in \cite{Smith1990-eq}.
However, this approach has not been applied to other networks as the monotonicity is not a general property of the DUE problems: 
indeed, the monotonicity does not hold in a simple network with only one origin-destination pair~\citep[][]{Kuwahara1990-vq,Mounce2007-ua}.
To address this issue, \cite{Satsukawa2019-lq} proved the stability in unidirectional networks~\citep[][]{Iryo2018-it} using a game-theoretic approach that does not rely on monotonicity.
This study showed that the DUE problems formulated as strategic games are weakly acyclic games based on the ordering property, and proved global convergence to a set of the equilibrium states and stochastic stability~\citep[][]{Young1993-bi}.
The applicability of this approach to more extended settings, such as general many-to-many networks with DTCs, is yet to be investigated.

\color{black}

\section{DTC game and its basic properties}\label{Sec:Preliminaries}

\subsection{Formulation of a DTC game}

\begin{figure}[t]
	\begin{center}
	\hspace{0mm}
    \includegraphics[width=130mm,clip]{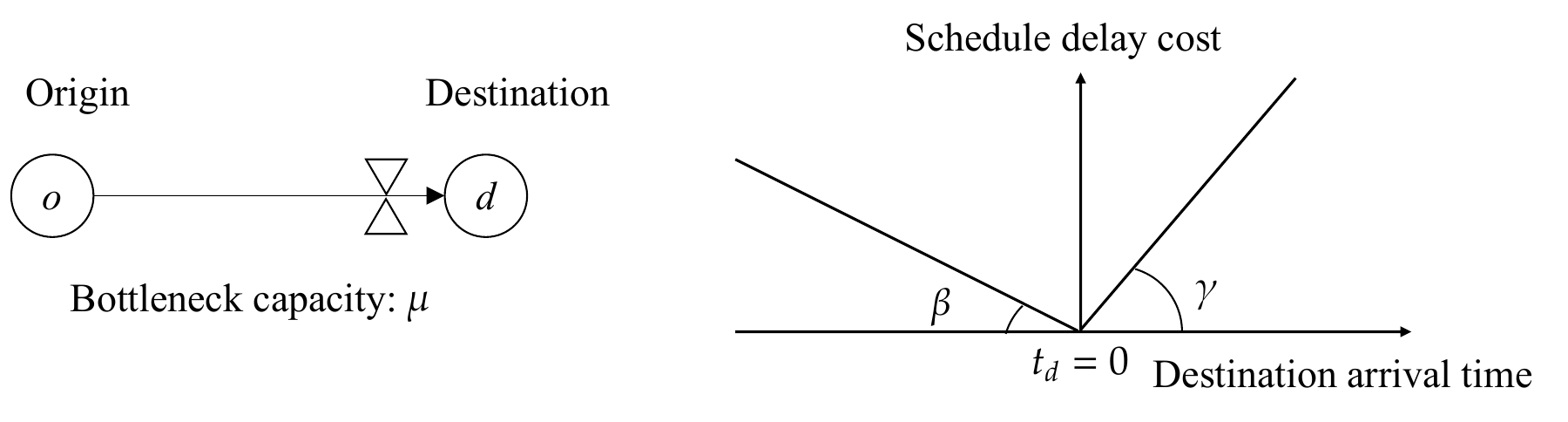}
	\end{center}
    \vspace{-4mm}
	\caption{Network and schedule delay cost function}
    \vspace{-2mm}
    \label{Fig:Network}
\end{figure}

We consider a road network with a single origin-destination pair connected by one link (Figure~\ref{Fig:Network}).
This link has a single bottleneck with a finite capacity of $\mu$ at the end.
A fixed number of homogeneous users travel along the link from the origin to the destination.
Each user is an atomic player in the DTC game.
The set of users is denoted by $\mathcal{P}\equiv\{1,\ldots, P\}$.
We define the size of the users to analyse the relationship between the DTC game and DTC problem of a conventional fluid model.
The size is denoted by $m$ $(0<m\leq 1)$.

Each user selects his/her origin departure time from a sufficiently long time period $[-S, S]$ \textcolor{black}{($S$ is a sufficiently large value)} so as to minimise the trip cost, i.e. each departure time is a strategy in the DTC game.
The time period is discretised with a sufficiently small value $\Delta s$, and users select their departure times from the following set, $[-S, -S+\Delta s, \ldots, 0, \ldots, S-\Delta s, S]$.
We set the value of $\Delta s$ so that users can precisely select departure times derived in the subsequent analysis of this study.\footnote{See \textbf{\ref{Sec:App-SettingDelta_t}} for the detailed setting of $\Delta s$ although the analysis includes variables derived in \textbf{Section~\ref{Sec:EpsilonNash}}.}
A strategy profile of the DTC game is a set of departure times of all users, referred to as a `time profile'.
A time profile is denoted by $\mathbf{s}\equiv\{s_{1},\ldots,s_{p},\ldots,s_{P}\}$, where $s_{p}$ represents the departure time of user $p\in\mathcal{P}$ in the time profile.
For any time profile $\mathbf{s}$, a set of the departure times of the users other than user $p$ is denoted by $\mathbf{s}_{-p}$.
With this notation, we sometimes represent the time profile $\mathbf{s}$ as $(s_{p},\mathbf{s}_{-p})$ to indicate the departure time of user $p$ clearly.

We assume that any user cannot select the departure times that the other users have already selected.
It follows that the set of available departure times for each user changes according to the time profile.
If the current time profile is $\mathbf{s}$, then the set of available departure times for user $p$ is denoted by $\mathcal{S}_{p}(\mathbf{s}_{-p})$.
We also denote by $\mathcal{S}$ the space of time profiles for the sake of convenience.

The trip cost, the (dis)utility of the DTC game, is assumed to be additively separable into the costs of free-flow travel, queueing delay and schedule delay.
The free-flow travel time is constant and assumed to be zero without loss of generality.
The queueing delay of a user is the difference between the departure time and the destination arrival time, which is uniquely determined for a given time profile based on a car-following model satisfying the first-in-first-out (FIFO) principle and causality~\citep[][]{Carey2003-qr}.
For consistency with the standard point queue model used in the fluid models, we assume that the destination arrival time is calculated as follows.
Let $o_{p}(\mathbf{s})$ denote the order of departure of the user $p\in\mathcal{P}$ from the origin when the time profile is $\mathbf{s}$.
We also denote by $d_{o_{p}(\mathbf{s})}(\mathbf{s})$ the destination arrival time of the user.
Then, the arrival time is determined by the following equation (hereinafter, we sometimes omit $\mathbf{s}$ from $o_{p}(\mathbf{s})$ and $d_{o_{p}(\mathbf{s})}(\mathbf{s})$ for simplicity):
\begin{align}
	&d_{o_{p}} = 
	\max\{ d_{o_{p}-1}+m/\mu, s_{p} \},\quad \text{where}\quad  d_{0} = -\infty.\label{Eq:QueueCost}
\end{align}
\noindent 
This means that the destination arrival time is determined such that the difference in the arrival time of its leading user (i.e. the headway) is maintained equal to or larger than $m/\mu$, which is the reciprocal of the bottleneck capacity considering the user size.
The first element in the maximum condition represents that the user must wait at the bottleneck to maintain the headway, i.e. the user travels in a congested situation;
the second element represents that the user travels in a free-flow situation.
Then, the destination arrival rate (i.e. the outflow rate from the bottleneck), which is the product of the user size and the reciprocal of the users' arrival interval, does not exceed the bottleneck capacity $\mu$; 
the queueing delay occurs when the arrival rate at the bottleneck exceeds $\mu$ in the same manner as in a point queue model.

The schedule delay cost is associated with the difference between the actual and desired arrival times at the destination.
Since all users are homogeneous, they have the same schedule delay cost function whose desired arrival time is $t_{d}=0$.
The schedule delay cost function is assumed to be piecewise linear, whose slopes for early and late arrivals are represented as rational numbers $\beta$ and $\gamma$, respectively (see the right-hand side of Figure~\ref{Fig:Network}).
Then, for a given destination arrival time $d$, the schedule delay cost function is defined as follows:
\begin{align}
    V(d) = \beta \max\{0 -d,0 \}+\gamma \max\{d,0 \}.
\end{align}

In summary, the trip cost of user $p\in\mathcal{P}$ for a given time profile $\mathbf{s}$ is expressed as follows:
\begin{align}
	C_{p}(\mathbf{s})=
	( d_{o_{p}}- s_{p} )
	+V(d_{o_{P}}).\label{Eq:TripCost}
\end{align}
\noindent Note that the value of time is assumed to be one without loss of generality.
We assume that $\beta<1$ according to the literature~\citep[e.g.][]{Hendrickson1981-bq}.\footnote{\cite{Hendrickson1981-bq} mentioned that `if a minute spent at the workplace is more onerous than a minute in queue, no stable equilibrium develops', and thus they set $\beta$ to a value less than the value of time.}

\subsection{Nash equilibrium and its non-existence issue}
One of the basic solution concepts for strategic games is pure-strategy Nash equilibrium.
The Nash equilibrium of the DTC game is defined as a traffic state (i.e. time profile) where no user can improve his/her trip cost by unilaterally changing the departure time.
Mathematically, a Nash equilibrium state $\mathbf{s}^{*}$ satisfies the following condition:
\begin{align}
	C_{p}(s^{*}_{p},\mathbf{s}^{*}_{-p}) = \min_{s\in\mathcal{S}_{p}(\mathbf{s}_{-p})}C_{p}(s,\mathbf{s}^{*}_{-p}),
	\quad \forall p\in\mathcal{P}.\label{Eq:DefiPureNash}
\end{align}
\noindent This definition is consistent with the equilibrium in the fluid models where no commuter has any incentive to change his/her departure time.
This implies that the trip costs of all homogeneous commuters are the same.

\begin{figure}[t]
	\begin{center}
	\hspace{0mm}
    \includegraphics[width=130mm,clip]{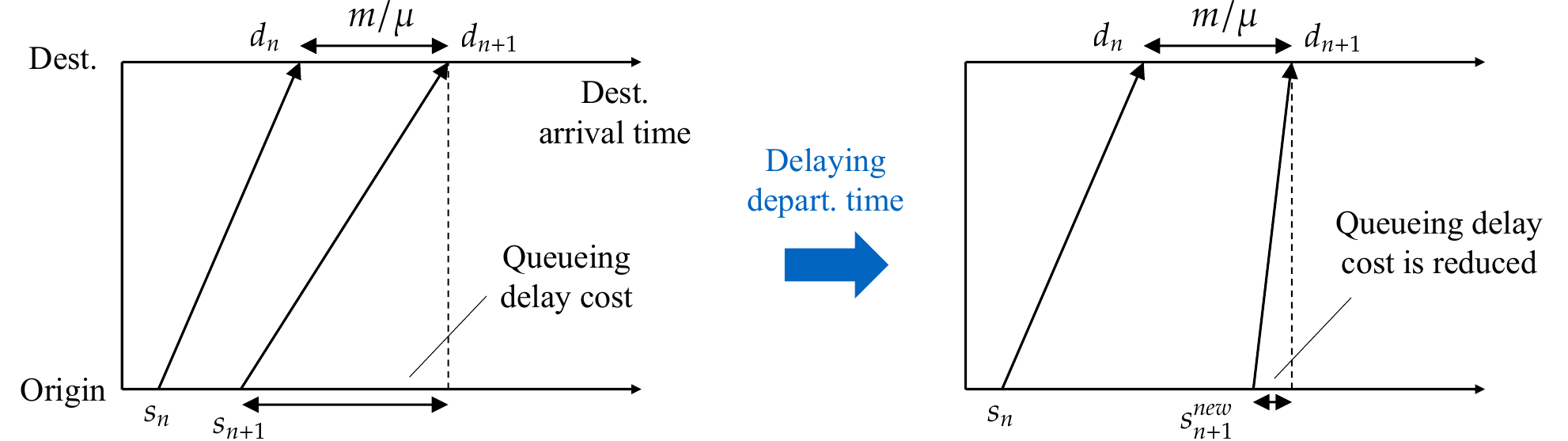}
	\end{center}
    \vspace{-4mm}
	\caption{Schematic of reducing queueing delay cost without changing the destination arrival time.}
    \vspace{-2mm}
    \label{Fig:QueueDelay}
\end{figure}

However, the pure-strategy Nash equilibrium does not always exist in DTC games because of the discretised definition of the users as pointed out in previous studies~\citep[][]{Otsubo2008-qj,Ziegelmeyer2008-hi,Ameli2022-zm}.
For example, suppose that users experience the same trip costs while experiencing positive queueing delay costs, as is the case with the equilibrium.
Eq.~\eqref{Eq:QueueCost} shows that their destination arrival times depend on the arrival times of the users in front of them.
Then, such a user can reduce the queueing delay cost without changing the destination arrival time, i.e. the schedule delay cost, by delaying his/her departure time (see Figure~\ref{Fig:QueueDelay}).
This suggests the non-existence of the pure-strategy Nash equilibrium.
Therefore, it is necessary to introduce another equilibrium concept to analyse the outcome of users' strategic behaviour. 


\section{Ensuring the existence: epsilon-Nash equilibrium}\label{Sec:EpsilonNash}
This section introduces the concept of \textit{epsilon-Nash equilibrium} to overcome the non-existence issue of the pure-strategy Nash equilibrium.
We first introduce the formal definition of the epsilon-Nash equilibrium.
We then prove the existence of this equilibrium in the DTC game.
Furthermore, we establish the correspondence between the epsilon-Nash equilibrium in the atomic model and the equilibrium in fluid models.

\subsection{Definition of epsilon-Nash equilibrium}
An epsilon-Nash equilibrium state is a strategy profile that approximately satisfies the Nash equilibrium condition with respect to an approximation parameter $\epsilon$.
Specifically, in an $\epsilon$-Nash equilibrium state, no user can improve his/her trip cost by more than $\epsilon$ by unilaterally changing his/her departure time.
This concept is formally defined as follows:
\begin{defi}
    Given a real non-negative parameter $\epsilon$, a time profile $\mathbf{s}^{*}$ is said to be $\epsilon$-Nash equilibrium of the DTC game if the following relationship holds:
    \begin{align}
    	C_{p}(\mathbf{s}^{*}) \leq C_{p}(s,\mathbf{s}^{*}_{-p}) + \epsilon,
    	\quad 
    	\forall s\in \mathcal{S}_{p}(\mathbf{s}^{*}_{-p}),\forall p\in\mathcal{P}.\label{Eq:DefiEpsilonNash}
    \end{align}
\end{defi}
\noindent The epsilon-Nash equilibrium corresponds to the pure-strategy Nash equilibrium if $\epsilon=0$.

The epsilon-Nash equilibrium is known as `bounded rational user equilibrium (BRUE)' in the transport research field.
In the BRUE, users do not have the incentive to change their strategies if their trip costs are within a certain value of the minimum trip cost~\citep[][]{Mahmassani1987-yp}.
This equilibrium concept has also been used to analyse the dynamic properties of the DTC problems with continuous traffic flows.
\cite{Guo2018-pu}, for example, developed day-to-day dynamics by extending the proportional swap system in \cite{Smith1984-ed} considering users' bounded rationality. 
\textcolor{black}{Although the convergence of the dynamics to the BRUE was not shown, they showed the convergence to the system optimal state under an appropriate pricing scheme.}
Similarly, \cite{Zhu2019-jv} analysed dynamical properties in a DTC problem with stochastic capacity, taking into account real-time information provision.
For a comprehensive review, see, for example, \cite{Di2016-eg}.



\subsection{Existence of the epsilon-Nash equilibrium}
Owing to the new equilibrium concept, we can ensure the existence of pure-strategy equilibrium in the DTC game.
Specifically, we prove that a time profile $\mathbf{s}^{e}$ satisfying the following conditions is epsilon-Nash equilibrium, inspired by the equilibrium in fluid models: the trip costs of all the users are the same, they arrive at the destination with the minimum headway $m/\mu$, and the queueing delay costs of the first and last users departing from the origin are zero.

\color{black}
We first derive the analytical solution to the identical trip cost in $\mathbf{s}^{e}$.
Let $t^{-}$ and $t^{+}$ denote the departure times of the first and last users, respectively.
Since their queueing delay costs are zero, their destination arrival times equal the departure times.
In addition, $P$ users arrive at the destination with the minimum headway $m/\mu$.
Combining these, we have 
\begin{align}
    t^{+} = t^{-} + \cfrac{m(P-1)}{\mu}.
\end{align}
We also have the following equation since the trip costs of the first and last users are the same:
\begin{align}
    -\beta t^{-} = \gamma \left(t^{-} + \cfrac{m(P-1)}{\mu}\right).
\end{align}
Solving these equations, we have the analytical solution to $t^{-}$, as follows:
\begin{align}
        &t^{-}=-\cfrac{m(P-1)}{\mu} \cfrac{\gamma}{\beta+\gamma}\label{Eq:DUE_ET}.
\end{align}
Then, the analytical solution to the trip cost $\rho$ is described in the following form:
\begin{align}
    \rho = \cfrac{m(P-1)}{\mu} \cfrac{\beta\gamma}{\beta+\gamma}.\label{Eq:DUE_Cost}
\end{align}

We next consider the departure time of each user $p\in\mathcal{P}$ whose order of departure is $o_{p}$.
Since the interval of destination arrival times is $m/\mu$, we have
\begin{align}
    d_{o_{p}} = t^{-} + \cfrac{m}{\mu}(o_{p}-1).
\end{align}
By substituting this into Eq.~\eqref{Eq:TripCost}, the trip cost of this user is described as follows:
\begin{align}
    C_{p} = t^{-} + \cfrac{m}{\mu}(o_{p}-1) - s^{e}_{p} + 
    \begin{cases}
        \beta \left\{ 0 - t^{-} - m(o_{p}-1)/\mu \right\}\quad 
        &\text{if}\quad t^{-} + \cfrac{m}{\mu}(o_{p}-1)\leq 0, \\
        \gamma \left\{ t^{-} +m(o_{p}-1)/\mu \right\}
        &\text{otherwise}.
    \end{cases}\label{Eq:App-Prf_As_1}
\end{align}
We also have the equation $-\beta t^{-} = C_{p}$ since all users have the same trip cost.
Then, substituting Eqs.~\eqref{Eq:DUE_ET} and \eqref{Eq:App-Prf_As_1} into this equation, we obtain the analytical solution to the departure time, $s_{p}^{e}$, as follows.
\begin{align}
    &s^{e}_{p}=t^{-}+
	\begin{cases}
		\cfrac{m(1-\beta)}{\mu}(o_{p}(\mathbf{s}^{e}) - 1)\quad 
		&\text{if}\quad o_{p}(\mathbf{s}^{e})\leq o^{cr}(\mathbf{s}^{e}),
		\\
		\cfrac{m(1+\gamma)}{\mu}(o_{p}(\mathbf{s}^{e})-1)-\cfrac{m \gamma (P-1)}{\mu}\quad 
		&\text{otherwise},
	\end{cases}\label{Eq:DUETimeProfile}
\end{align}
where $o^{cr}(\mathbf{s}^{e})$ is the number of users who arrive at the destination no later than the desired arrival time for the given time profile $\mathbf{s}^{e}$ and is represented by
\begin{align}
    o^{cr}(\mathbf{s}^{e}) = \lfloor \gamma(P-1)/(\beta+\gamma) \rfloor +1.
\end{align}
\color{black}

\noindent The first case in Eq.~\eqref{Eq:DUETimeProfile} describes the departure times of early arrival users whose destination arrival times are equal to or not later than zero, while the second case describes those of late arrival users whose destination arrival times are later than zero.
This suggests that the departure time interval between users is equal to or less than $m(1+\gamma)/\mu$: 
the interval between early arrival users is $m(1-\beta)/\mu$, between late arrival users is $m(1+\gamma)/\mu$ and the interval between the latest early arrival and the earliest late arrival users is equal to or less than $m(1+\gamma)/\mu$.

\begin{figure}[t]
	\begin{center}
	\hspace{0mm}
    \includegraphics[width=130mm,clip]{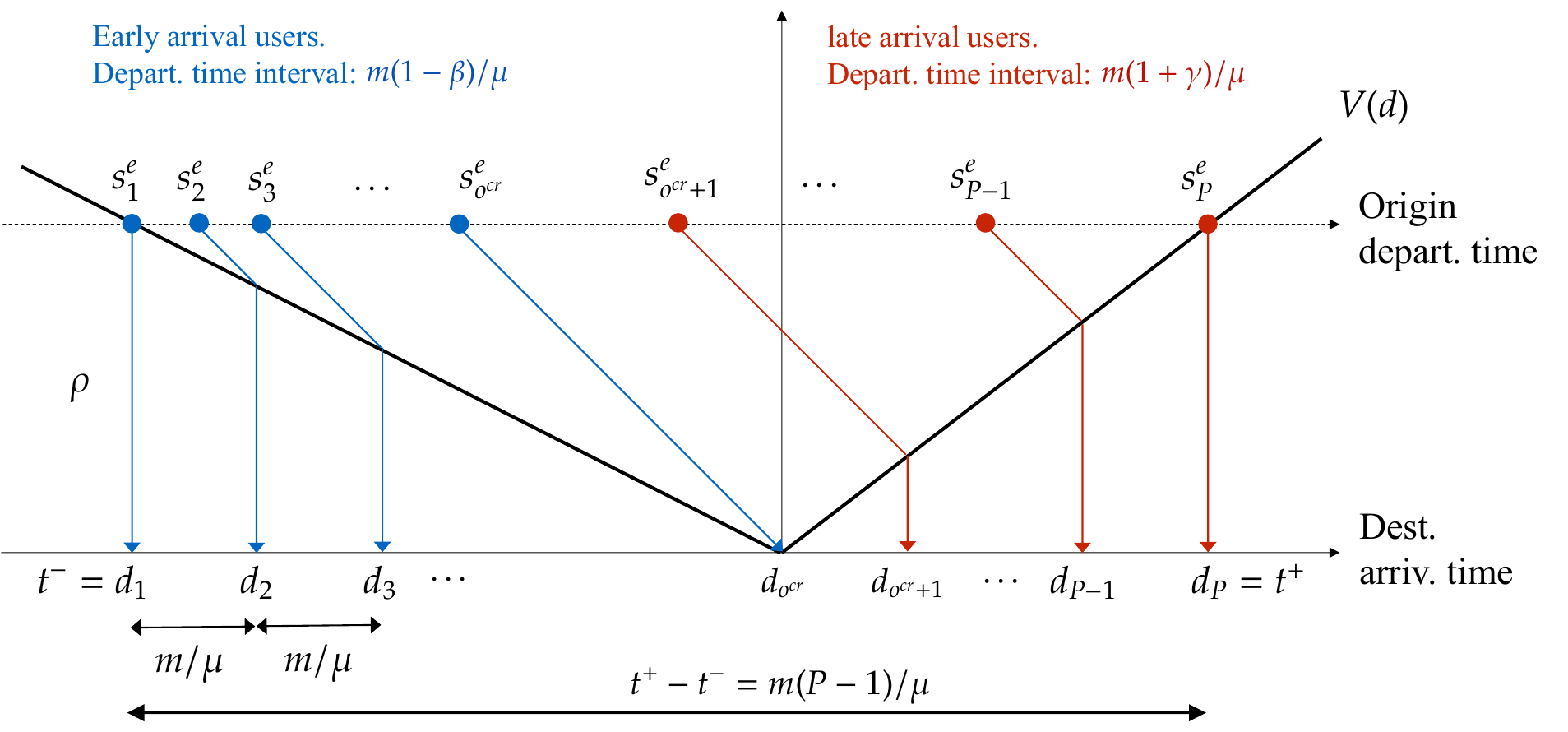}
	\end{center}
    \vspace{-4mm}
	\caption{Description of the pattern of the departure and destination arrival times in $\mathbf{s}^{e}$; the interval between destination arrival times is $m/\mu$ and users experience corresponding queueing and schedule delay costs, of which the sum is $\rho$.}
    \vspace{-2mm}
    \label{Fig:EpsilonNash}
\end{figure}

The important point is that the departure time interval becomes the upper limit for the improvement in the trip cost of a user by unilaterally changing the departure time.
In $\mathbf{s}^{e}$, the interval of the destination arrival times between consecutive users is $m/\mu$.
This means that users travel in congested situations, and their destination arrival times become dependent on those of the users in front of them.
It follows that a user can reduce his/her queueing delay cost without changing the schedule delay cost by delaying the departure time, i.e. by increasing the interval of departure times between the user and its leader user, as mentioned in the discussion on the non-existence issue in \textbf{Section~\ref{Sec:Preliminaries}}.
Therefore, the departure time interval becomes an important factor in improving the users' trip costs.
Note that the formal analysis in the upper limit has to deal with complicated situations, such as when a user's departure order varies due to the change in the departure time, but the logic of the improvement in such situations is basically the same.


Based on the above discussion, we can prove that the time profile $\mathbf{s}^{e}$ becomes the epsilon-Nash equilibrium if the parameter $\epsilon$ is sufficiently large, as follows:
\begin{theo}\label{Theo:ExistenceNash}
	Suppose that $\epsilon$ in Eq.~\eqref{Eq:DefiEpsilonNash} is equal to or higher than $m(1+\gamma)/\mu$.
	Then, the time profile $\mathbf{s}^{e}$, in which the departure time of each user $p\in\mathcal{P}$ is defined by Eq.~\eqref{Eq:DUETimeProfile}, is the $\epsilon$-Nash equilibrium of the DTC game.
\end{theo}
\color{black}
\begin{prf}
	 The outline of the proof is shown in the above paragraph. See \textbf{\ref{App:Prf_Improvement}} for the detailed analysis.\qed
\end{prf}
\color{black}

\noindent This theorem ensures the existence of pure-strategy equilibrium in the DTC game.
Hereinafter, we consider the case where $\epsilon = m(1+\gamma)/\mu$ for the given parameter $m$ in order to avoid some unnecessary complications for the analysis.

It is worth noting that the epsilon-Nash equilibrium in a DTC game would not be unique: 
there could exist other equilibrated time profiles where users experience different trip costs since the parameter $\epsilon$ allows users to select the non-best departure times in equilibrium, unlike the equilibrium in DTC problems with a fluid approach.
Nevertheless, it is valuable to focus on this specific epsilon-Nash equilibrium $\mathbf{s}^{e}$ because this equilibrium becomes a key for investigating the relationship with the equilibrium in fluid models, as will be shown in the next section.


\subsection{\textcolor{black}{Relationship to the equilibrium in fluid models}}

\subsubsection{Revisit to the equilibrium in Vickrey's continuous bottleneck model}
We first revisit the equilibrium in the continuous bottleneck model of \cite{Vickrey1969-rg}.
Consider a DTC problem in a single-road network with the same bottleneck capacity $\mu$.
Commuters are treated as a continuum, and the total mass is denoted by $Q$.
All commuters are homogeneous: they have the same desired arrival time, value of time and schedule delay cost function as with the users in the DTC game.
The queueing delay at the bottleneck is described by a standard point queue model in accordance with the FIFO principle and causality.
This cost for commuters departing from the origin at time $t$ is mathematically expressed as follows~\citep[][]{Hendrickson1981-bq,Arnott1990-oa}:
\begin{align}
    &T(t) = \cfrac{D(t)}{\mu},\quad 
    \text{where}\quad 
    D(t) = \int_{{\color{black} \hat{t}}}^{t}r(u)\mathrm{d}u 
    - \mu(t-\hat{t}).
\end{align}
\noindent Here, $r(u)$ denotes the departure flow rate function and \textcolor{black}{$\hat{t}$ denotes the most recent time at which there was no queue}.

The equilibrium is defined as the state where no commuter can reduce their trip cost by changing their departure time, as mentioned in \textbf{Section~\ref{Sec:Preliminaries}}.
We denote by $\rho^{f}$ the equilibrium trip cost and $t^{f}$ the departure time of the first commuter in the equilibrium.
These variables are analytically given by
\begin{align}
    &\rho^{f} = \cfrac{Q}{\mu}\cfrac{\beta\gamma}{\beta+\gamma},\quad t^{f} = -\cfrac{Q}{\mu}\cfrac{\gamma}{\beta+\gamma}.\label{Eq:fluidDUE_Cost}
\end{align}
Moreover, the departure flow rate for early and late arrivals, which are denoted by $r^{e}$ and $r^{l}$ respectively, are piecewise constant and given by 
\begin{align}
    r^{e} = \cfrac{\mu}{(1-\beta)},\quad r^{l} = \cfrac{\mu}{1+\gamma}.
\end{align}

\subsubsection{Correspondence between the epsilon-Nash equilibrium and the equilibrium in the DTC problem}
By comparing Eqs.~\eqref{Eq:DUE_ET} and \eqref{Eq:DUE_Cost} to Eq.~\eqref{Eq:fluidDUE_Cost}, we have the following proposition whose proof is straightforward and thus omitted:
\begin{prop}
    Consider the DTC game with $P$ users and the DTC problem with $Q$ total mass.
    Suppose that the number of users and the total mass satisfy the following relationship: $m(P-1)=Q$.
    Then, the identical trip costs of all users in the epsilon-Nash equilibrium $\mathbf{s}^{e}$, $\rho$, become the same as the equilibrium trip cost in the equilibrium, $\rho^{f}$.
    In addition, the departure times of the first and last users (i.e. the rush-hour period) are exactly the same.
\end{prop}
\noindent This proposition claims that the trip cost $\rho$ in the epsilon-Nash equilibrium $\mathbf{s}^{e}$ becomes the same as $\rho^{f}$ by appropriately adjusting the number of users and their size.
Note that the departure flow rates in $\mathbf{s}^{e}$ also correspond to those in the equilibrium in the DTC problem.
Since the departure time interval for early arrival users is $m(1-\beta)/\mu$ when the user size is $m$, the departure flow rate is give by
\begin{align}
    \cfrac{1}{m(1-\beta)/\mu}\cdot m = \cfrac{\mu}{1-\beta}.
\end{align}
The departure flow rate for late arrival users is calculated as $\mu/(1+\gamma)$ similarly.
These results confirm that the departure flow rates in $\mathbf{s}^{e}$ and the equilibrium in the DTC problem are the same.

Moreover, this proposition suggests that the epsilon-Nash equilibrium $\mathbf{s}^{e}$ asymptotically approaches the equilibrium in the DTC problem as the user size $m$ approaches zero while the relationship $m(P-1)=Q$ is satisfied.
As $m\rightarrow 0$, the value of $\epsilon = m(1+\gamma)/\mu$ ensuring the existence of the equilibrium $\mathbf{s}^{e}$ approaches zero.
This means that as the size approaches zero while the population size accordingly increases, $\mathbf{s}^{e}$ approaches the equilibrium in the DTC problem where every user selects their \textit{truly best} departure time at which they cannot improve their trip cost.
Therefore, the epsilon-Nash equilibrium enables us to establish not only the existence of pure-strategy equilibrium but also its correspondence to the equilibrium of the fluid model, i.e. the correspondence between the atomic and fluid models.\footnote{We can prove that the epsilon-Nash equilibrium states other than $\mathbf{s}^{e}$ also approach the equilibrium in the DTC problem as $m\rightarrow 0$ although the proof is complicated and thus omitted.}

\section{Global convergence of the DTC game}\label{Sec:WAG}

\color{black}

This section analyses the global convergence to the equilibrium using the following two concepts: weakly acyclic games and the ordering property.
First, we introduce the concept of a weakly acyclic game, which is a class of games that satisfies the necessary conditions for global convergence.
We next prove that the DTC game is a weakly acyclic game by showing the existence of a better response path that converges to the epsilon-Nash equilibrium based on the ordering property.
By analysing the convergent better response path, we clarify the user behaviour necessary to achieve the global convergence, i.e. the convergence mechanisms.
Finally, based on the clarified convergence mechanisms, we develop the evolutionary dynamics under which global convergence is guaranteed.
We also discuss the reason for the stability and instability of the existing evolutionary dynamics studied in the literature.

\color{black}


Hereinafter, we consider a situation where a DTC game is played repeatedly on a day-to-day basis.
Let $\tau \in \mathbb{N}$ denote the time step index (i.e. day).
We denote by $s_{p}^{\tau}$ the departure time of a user $p\in\mathcal{P}$ on a day $\tau$.
In such a repeated game, some users are randomly selected on each day.
They change their departure times according to a behavioural rule while referring to the trip costs at other departure times.
The other users continue to choose their departure time.
The resulting (day-to-day) evolution of the traffic state is called evolutionary dynamics.

\subsection{Weakly acyclic games: definition and its importance}
We first introduce the concept of an \textit{improvement path} that characterises a class of weakly acyclic games.
Consider a sequence of strategy (time) profiles $\mathbf{s}^{1},\ldots,\mathbf{s}^{\tau},\ldots,\mathbf{s}^{L}$.
Following \cite{Monderer1996-ov}, a \textit{path} is a sequence such that, for each day $\tau>1$, there is exactly one user $p\in\mathcal{P}$ who changes the strategy, i.e. $\mathbf{s}^{\tau} = (s_{p},\mathbf{s}_{-p}^{\tau-1})$ and $s_{p}\neq s_{p}^{\tau-1}$.
Then, a path is called an improvement path if the change in strategy improves the utility (trip cost) forecasted by the user.
Mathematically, each successive pair $(\mathbf{s}^{\tau-1},\mathbf{s}^{\tau})$, $\forall \tau>1$ in an improvement path satisfies the following relationship:
\begin{align}
    C_{p}(\mathbf{s}^{\tau-1}) > \hat{C}_{p}(s^{\tau}_{p}\mid \mathbf{s}^{\tau-1}),\label{Eq:ImprovementPath}
\end{align}
where $\hat{C}_{p}(s' \mid \mathbf{s})$ denotes the trip cost for the departure time $s'$ forecasted by user $p$ when the current time profile is $\mathbf{s}$.
We refer to this trip cost as a forecasted cost (utility).

An improvement path is specified according to the definition of the forecasted utility.
For example, \cite{Marden2009-bs} considered a \textit{better improvement path} by assuming that the forecasted utility for a strategy is equal to the utility which a user actually experiences after the change to the strategy;
that is, Eq.~\eqref{Eq:ImprovementPath} in the better improvement path is described as follows:
\begin{align}
    C_{p}(\mathbf{s}^{\tau-1}) > C_{p}(s^{\tau}_{p},\mathbf{s}_{-p}^{\tau-1}).
\end{align}
\cite{Young1993-bi} also considered a \textit{best improvement path} in which one user changes the strategy so as to maximise the future utility at each time step.

We then introduce the definition of weakly acyclic games based on a specific improvement path, as follows~\citep[][]{Young1993-bi,Fabrikant2013-cp}:
\begin{defi}
    A game is a weakly acyclic game (under the specified improvement path) if (i) for any strategy profile $\mathbf{s}\in\mathcal{S}$ that is not in equilibrium, there exists an improvement path starting at $\mathbf{s}$ and ending at a Nash equilibrium state, and (ii) there is no improvement path from the Nash equilibrium state, i.e. that equilibrium is a stationary point.
\end{defi}
\noindent 
By definition, weakly acyclic games are positioned as game classes that satisfy the necessary condition for a collective state to reach Nash equilibrium through users' rational and distributed behaviour, i.e. for the convergence of evolutionary dynamics~\citep[][]{Fabrikant2013-cp}.
The behaviour of `natural' evolutionary dynamics, such as satisfying positive correlation property~\citep[][]{Sandholm2010-ht}, can be described by some improvement path.
This implies that the evolutionary dynamics could converge to equilibrium states of a strategic game if the game is a weakly acyclic game under the corresponding improvement path.
It is therefore important for the analysis of dynamic properties in strategic games to check whether a game is a weakly acyclic game or not.\footnote{Indeed, \cite{Sandholm2001-yv,Sandholm2002-wg} proved the convergence and stability in congestion games by linking them to potential games, which are special cases of weakly acyclic games.}
Conversely, if a game is not a weakly acyclic game, we cannot guarantee global convergence and stability.
There exists at least one state that has no improvement path leading to an equilibrium state (i.e. equilibrium cannot be reached from that state).

This definition of a weakly acyclic game is slightly more general than that usually used in game theory.
In this research field, the forecasted utility is generally defined as the utility that a user actually experiences after the strategy change, as shown in the examples of improvement paths.
The reason why we introduce the general definition is that it would be reasonable to assume that users may not know such a future utility in the field of traffic assignment, where the population is so large that the players do not fully know the structure of the game.
Furthermore, this definition enables us to prove the weakly acyclicity of the DTC game, while we cannot prove it under the original definition.


\subsection{DTC game is a weakly acyclic game with a better response path}


\subsubsection{Weakly acyclicity of DTC games}
We first define the forecasted cost for a departure time when no user is departing, in order to define an improvement path in this study.
For this issue, we assume that the forecasted cost is expressed by appropriate interpolation based on the current trip costs and schedule delay costs of the users.

Consider a situation where a user $p\in\mathcal{P}$ changes his/her departure time to a time $s'$ between the departure times of the two consecutive users $a,b\in\mathcal{P}$, \textcolor{black}{i.e. $s'\in(s_{a},s_{b})$.}
If they travel in congested situations (i.e. the interval between their destination arrival times is $m/\mu$), the user forecasts that $s'$ is included in a congested time interval, i.e. a user departing at that time will travel in a congested situation.
Then, the forecasted cost for $s'$ can be determined by the linear interpolation from the departure times and trip costs of the two congested users, as follows: 
\begin{align}
    \hat{C}_{p}(s' \mid \mathbf{s}) = 
    C_{a} + \cfrac{C_{b} - C_{a}}{s_{b} - s_{a}}(s'-s_{a}).\label{Eq:DefiBetterProto}
\end{align}
Meanwhile, if the two users travel in free-flow situations (i.e. the interval between their destination arrival times is larger than $m/\mu$), the time interval $(s_{a},s_{b})$ may have a congested time interval and a free-flow time interval, with the time $d_{o_{a}}$ (i.e. queue dissipation time) as the boundary.
Then, if $s'\in (s_{a}, d_{o_{a}} ]$ i.e. included in a congested time interval, the forecasted cost is determined from Eq.~\eqref{Eq:DefiBetterProto} with replacing $(s_{b},C_{b})$ by $(d_{o_{a}}, V(d_{o_{a}}) )$; 
if $s'\in (d_{o_{a}}, s_{b})$, i.e. included in a free-flow time interval, the forecasted is described as the corresponding schedule delay cost $V(s')$.
Note that when $s'$ is earlier than the departure time of the first user departing from the origin or later than the destination arrival time of the last user, the forecasted cost is also defined as $V(s')$.

This definition is consistent with standard day-to-day dynamics used in fluid models, under which each user changes his/her departure time with reference to the actual costs realised in the current traffic state.
We then define the following better response path, on which users change their departure times so as to strictly improve the utility based on the forecasted costs: 
\begin{defi}\label{Defi:BRP}
    An improvement path is a better response path if users change their departure times according to Eq.~\eqref{Eq:ImprovementPath} with the forecasted costs.
\end{defi}

\color{black}
We now prove that the DTC game is a weakly acyclic game based on the following \textit{ordering property}~\citep[][]{Iryo2011-bk,Satsukawa2019-lq}: 
\begin{defi}
The ordering property is the existence of the appropriate assignment order such that equilibrium is achieved by changing users' strategies to the equilibrium strategies according to the order.
\end{defi}

\noindent Referring to the ordering property, we prove the existence of a convergent better response path where users change their strategies in order, as follows (we call this `ordered path'): 
users' departure times are sequentially changed to equilibrium departure times in $\mathbf{s}^{e}$, starting from the earliest equilibrium ones;
in addition, such changes are made within the time range later than the departure times of the equilibrated users who have already selected the equilibrium strategies.
On such an ordered path, once users change their departure times to the equilibrium ones, their trip costs become the equilibrium one, $\rho$.
Furthermore, their trip costs are not changed by the change in the departure times of non-equilibrated users, who depart later than the equilibrated users, due to the FIFO principle and causality of the bottleneck model: a user's link travel time is not affected by the behaviour of users who depart later.
Therefore, the sequential change will achieve the epsilon-Nash equilibrium in which all users experience the same trip costs.

Based on this idea, we prove the following theorem:

\color{black}
\begin{theo}\label{Theo:WAG}
    A DTC game is a weakly acyclic game under the better response path of \textbf{Definition~\ref{Defi:BRP}} in the following sense: (i) from any time profile $\mathbf{s}\in\mathcal{S}\setminus\{\mathbf{s}^{e}\}$, there exists a better response path starting at the time profile and ending at the epsilon-Nash equilibrium $\mathbf{s}^{e}$, and (ii) there is no better response path from the epsilon-Nash equilibrium state.
\end{theo}
\begin{prf}
The outline of the proof is described as follows.
We first establish the following lemma:

\begin{lemm}\label{Lemm:TheoWAG-1}
Consider a time profile $\mathbf{s}$ where the first $n$ users $(1\leq n < P)$ departing from the origin select the corresponding equilibrium departure times in Eq.~\eqref{Eq:DUETimeProfile}, i.e. $s_{p} = s_{p}^{e}$, if $o_{p}(\mathbf{s})\leq n$.
Then, there exists a better response path from $\mathbf{s}$ to a time profile $\mathbf{s}'$ where the first to $n+1$st users select the equilibrium departure times. 
Moreover, during the sequence of better responses, (a) the first $n$ equilibrated users do not change their departure times, and (b) the other (non-equilibrated) users do not change their departure times to earlier times than the equilibrated users.

\end{lemm}

\begin{figure}[t]
	\begin{center}
	\hspace{0mm}
    \includegraphics[width=130mm,clip]{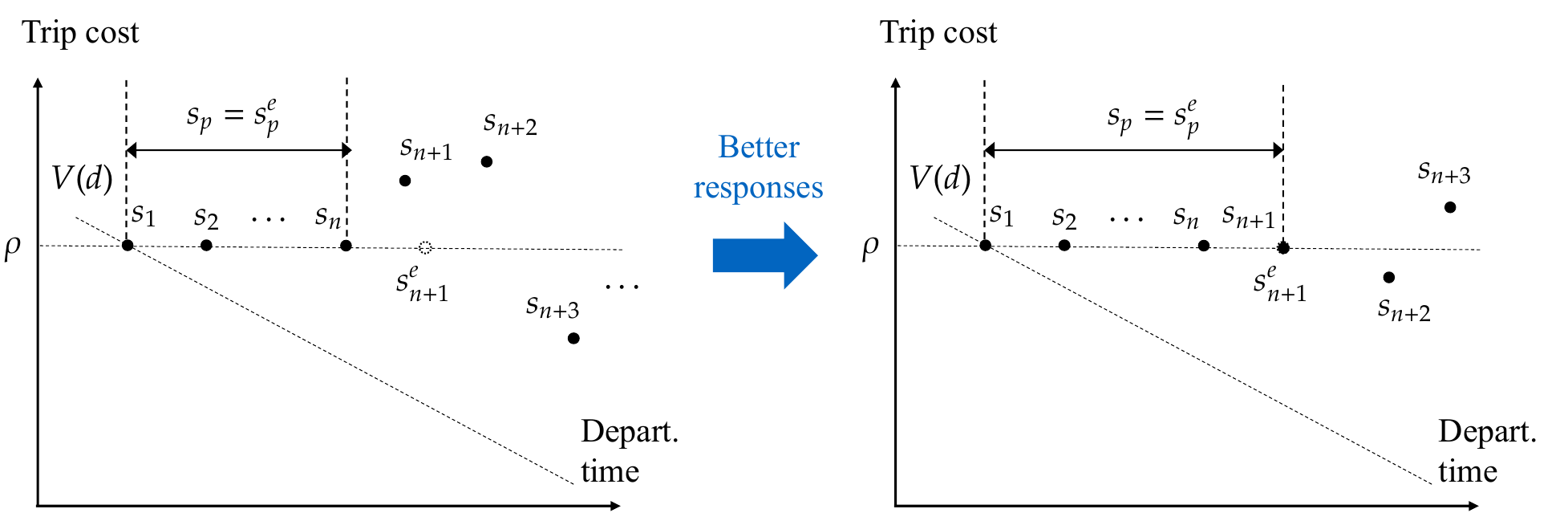}
	\end{center}
    \vspace{-4mm}
	\caption{Schematics of a better response of $n+1$st user to $s_{n+1}^{e}$; each dot represents the relationship between a departure time and corresponding trip cost of each user. 
    }
    \vspace{-2mm}
    \label{Fig:BRP}
\end{figure}


This lemma shows that when the first user selects the equilibrium departure time $t^{-}$, subsequent users can be fixed to the equilibrium departure times in order through better responses (see Figure~\ref{Fig:BRP}).
In addition, this lemma guarantees that the trip costs of such equilibrated users are also fixed to the equilibrium trip cost $\rho$:
since a non-equilibrated user does not select a departure time earlier than those of the equilibrated users, their destination arrival times (i.e. queueing and schedule delay costs) are not changed by the DTCs of non-equilibrated users.
This means that the equilibrium $\mathbf{s}^{e}$, where all users experience the same trip cost $\rho$, can be realised through better responses where users change their departure times so as to experience the same trip cost as that of the first user.
This ensures the existence of a better response path to $\mathbf{s}^{e}$ from a time profile $\mathbf{s}$ where the first user selects $t^{-}$.

Next, we prove the existence of a better response path from an arbitrary initial time profile to such a time profile, where the first user selects the equilibrium departure time $t^{-}$.
To this end, we establish the following lemma:
\begin{lemm}\label{Lemm:TheoWAG-2}
Consider a time profile $\mathbf{s}$ where the departure time of the first user $s_{1}$ is not $t^{-}$.
Then, there exists a better response path to a time profile $\mathbf{s}'$ where the departure time of the first user $s'_{1}$ satisfies the following relationship: (i) if $s_{1}<t^{-}$, then $s_{1}<s'_{1}$ and (ii) if $t^{-}<s_{1}$, then $s'_{1} = t^{-}$.
\end{lemm}

\noindent This lemma shows that the departure time of the first user $s_{1}$ can be adjusted to $t^{-}$ through better responses.
Specifically, if $s_{1}$ is earlier than $t^{-}$ in the current time profile $\mathbf{s}$, then the first user can conduct a better response to another departure time later than $s_{1}$, which realises the new time profile $\mathbf{s}'$ where the departure time of the first user $s_{1}'$ is later than $s_{1}$.
Meanwhile, if $s_{1}$ is later than $t^{-}$, then some user can conduct a better response to $t^{-}$, which means that the departure time of the first user becomes the same as $t^{-}$ in the new time profile $\mathbf{s}'$.

\textcolor{black}{
Consequently, \textbf{Lemma~\ref{Lemm:TheoWAG-1}} and \textbf{Lemma~\ref{Lemm:TheoWAG-2}} ensure that the epsilon-Nash equilibrium $\mathbf{s}^{e}$ can be achieved by sequential changing the departure times of non-equilibrated users to the equilibrium departure times in order, i.e. the existence of an ordered path.
This means the existence of a better response path to $\mathbf{s}^{e}$ from an arbitrary initial time profile in the DTC game, other than $\mathbf{s}^{e}$.
}


Furthermore, it is obvious that the forecasted cost for any departure time is equal to or higher than $\rho$ in $\mathbf{s}^{e}$ for the following reasons.
Since the interval between the destination arrival times of any two consecutive users is $m/\mu$ and all users experience the same trip cost $\rho$, the forecasted cost for a departure time in $(t^{-}, t^{+})$ is calculated from Eq.~\eqref{Eq:DefiBetterProto} as $\rho$.
In addition, the forecasted cost for another departure time becomes equal to the corresponding schedule delay cost, which is obviously higher than $\rho$.
Therefore, $\mathbf{s}^{e}$ is the unique stationary point of the better responses.
These prove the theorem (see \textbf{\ref{App:Prf_WAG}} for the proofs of the lemmas).
\qed

\end{prf}

\subsubsection{Mechanisms underlying the convergence of the better responses}
\textcolor{black}{
\textbf{Theorem~\ref{Theo:WAG}} enables us to gain insight into the convergence mechanisms from the ordered path, which is a better response path that converges to equilibrium.
This theorem states that users' better responses converge to the epsilon-Nash equilibrium $\mathbf{s}^{e}$ from any initial time profile by following the appropriate assignment order.
Thus, by analysing the convergent better response path, we can find important behavioural rules that ensure global convergence to the equilibrium. 
Specifically, the lemmas in the proof tell us that the following behavioural rules are required to achieve global convergence by tracing the ordered path.
}



\textbf{Lemma~\ref{Lemm:TheoWAG-1}} indicates that the departure times of users are sequentially fixed to the equilibrium departure times so as to experience the same (equilibrium) trip cost as that of the first user $\rho$.
There are two types of mechanisms behind this fixation, (a) and (b), as shown in the lemma.
Mechanism (a) prevents the equilibrated users from further changing their departure times.
Non-equilibrated users, who depart later than the equilibrated users, may experience trip costs lower than $\rho$;
then, equilibrated users would have an incentive to change their departure times to neighbourhoods of non-equilibrated users.
In other words, the ordered path has a property for maintaining a traffic state close to the equilibrium by preventing the equilibrated users from delaying their departure times.

Mechanism (b) limits the departure times available for non-equilibrated users to departure times later than those of the equilibrated users.
This can prevent the following `ripple effect'.
A better response of a non-equilibrated user to a departure time earlier than equilibrated users (i.e. overtaking) will delay their destination arrival times.
This increases their queueing delay costs and the trip costs become higher than the equilibrium cost $\rho$, which makes them non-equilibrated users.
Then, the non-equilibrated users will have the incentive to change their departure times, and their better responses to earlier departure times than (remaining) equilibrated users will move the traffic state away from the equilibrium.
In other words, the ordered path has a property for preventing such a perturbation due to the better responses with overtaking.

\textbf{Lemma~\ref{Lemm:TheoWAG-2}} shows that the departure time of the first user can be asymptotically adjusted to the equilibrium departure time $t^{-}$ through better responses.
Consider the difference in the departure time of the first user between the current and equilibrium time profiles, $\Delta t = s_{1} - t^{-}$.
Then,  when $s_{1}$ is earlier than $t^{-}$, a better response of the first user can achieve a new time profile $\mathbf{s}'$ where the departure time of the first user $s'_{1}$ is later than $s_{1}$, which means that $\Delta t$ approaches zero;
once $s_{1}$ becomes later than $t^{-}$, i.e. $\Delta t < 0$, an appropriate better response can lead $s'_{1}$ to $t^{-}$, which means that $\Delta t = 0$.
This implies that in the convergent better response path, the upper and lower bounds of the difference in the departure time of the first user can approach zero.


These results suggest the importance of the fixation of the equilibrated users and the asymptotic adjustment of the departure time of the first user for \textcolor{black}{tracing the ordered path}.
Based on this, the next section develops evolutionary dynamics that globally converges to the equilibrium.

\subsection{Convergence and stability properties in DTC games}

\subsubsection{Basic better response dynamics}
This section proceeds with the analysis based on the following \textit{better response dynamics}.
On each day $\tau$, one user $p\in\mathcal{P}$ is selected randomly with an equal probability (i.e. $1/P$).
The selected user randomly selects a departure time $s'$ and calculates the forecasted cost $\hat{C}_{p}(s' \mid \mathbf{s}^{\tau})$.
Then, if $C(\mathbf{s}^{\tau}) > \hat{C}(s' \mid \mathbf{s}^{\tau})$ holds, the user changes to $s'$; otherwise, the user does not change the departure time.
In other words, user $p$ changes the departure time by selecting one from the following set of better responses:
\begin{align}
    D_{p}(\mathbf{s}^{\tau}):=
    \left\{ s^{*} \mid s^{*}\in\mathcal{S}_{p}(\mathbf{s}_{-p}) \ \text{s.t.}\ C_{p}(\mathbf{s}^{\tau}) > \hat{C}(s' \mid \mathbf{s}^{\tau}) \right\}.
\end{align}

The better response dynamics can take the ordered path leading to the equilibrium.
Therefore, we can theoretically establish the \textit{almost sure convergence} properties of the dynamics as a direct consequence of \textbf{Theorem~\ref{Theo:WAG}}, as in \cite{Satsukawa2019-lq,Satsukawa2022-zn}.
However, such naive evolutionary dynamics is assumed to have difficulties converging to equilibrium.
This is because the ordered path seems to be rarely selected under the evolutionary dynamics since this does not have the mechanisms for the convergence pointed out in the previous section.

\subsubsection{Better response dynamics with the fixation and asymptotic adjustment}

Considering the convergence mechanisms, we formulate the following variation of the better response dynamics that repeat the fixation and adjustment processes of departure times of users:

\begin{description}
    \item[Step 0] 
    Let also $\overline{s}$ and $\underline{s}$ denote the maximum and minimum possible departure time for the first user, respectively: the departure time of the first user is fixed if he/she changes to the time within the range of $(\underline{s},\overline{s})$.
    These variables are initially set sufficiently large and small, e.g. $\overline{s} = S$ and $\underline{s} = 0$.
    Set initially $\tau = 0$.

    \item[Step 1] 
    Consider the current time profile $\mathbf{s}^{\tau}$.
    Let $s_{1}$ and $C^{r}$ denote the departure time and trip cost of the first user, respectively.
    Suppose that the first $n$ users experience the same trip cost $C^{r}$.
    The set of the other users is denoted by $\overline{\mathcal{P}}$.

    \item[Step 2] On each day $\tau$, one user is randomly selected from the set $\overline{\mathcal{P}}$ with equal probability (i.e. $1/|\overline{\mathcal{P}}|$).
    The selected user is denoted by $p\in\overline{\mathcal{P}}$.
    User $p$ randomly selects a new departure time from the set $D_{p}(\mathbf{s}^{\tau})$ not to depart earlier than the first $n$ users.

    \item[Step 3] If the $n+1$st user $p'$ experiences the trip cost $C^{r}$, then $\overline{\mathcal{P}}:=\overline{\mathcal{P}}\setminus\{p'\}$.
    On the next day $\tau:=\tau+1$, one of the users in $\overline{\mathcal{P}}$ again conducts a better response as explained in Step 2.

    \item[Step 4] If no user in $\overline{\mathcal{P}}$ cannot conduct such a better response so as to experience the trip cost $C^{r}$, the trip costs of users are checked.
    
    If the trip costs of the users in $\overline{\mathcal{P}}$ are higher than $C^{r}$, then $\underline{s}$ is updated as follows: $\underline{s}:=s_{1}$;
    If the trip costs of the users in $\overline{\mathcal{P}}$ are lower than $C^{r}$ or the last user experiences a positive queueing delay (i.e. the user can conduct a better response), then $\overline{s}$ is updated as follows: $\overline{s}:=s_{1}$.

    After the update, $\overline{\mathcal{P}}:=\mathcal{P}$, i.e. all users become able to conduct better responses.
    If the departure time of the first user changes to a time in $(\underline{s},\overline{s})$, update $C^{r}$ to the schedule delay cost of the first user.
    Then go back to Step 1 and repeat better responses.
\end{description}

\noindent 
This dynamics searches for the equilibrium time profile $\mathbf{s}^{e}$ without knowing it in advance by repeating the following processes: the departure times of users are fixed so as to experience the same trip cost as that of the first user, assuming that the user selects the equilibrium departure time; 
if not, the departure time of the first user is appropriately adjusted based on the resulting time profile.

Step 0 corresponds to the initialisation which prepares for the repetition of the dynamics while adjusting the departure time of the first user.
Step 1 corresponds to the preparation for fixing the departure times of users.
In this step, the current first user becomes a reference user, who is assumed to have selected the equilibrium departure time.
Subsequent users sequentially adjust their departure times to experience the same trip cost as $C^{r}$ of the reference user in Steps 2 and 3.
This implies that this better response will converge to $\mathbf{s}^{e}$ if the first user indeed selects the equilibrium departure time $t^{-}$, as mentioned in \textbf{Lemma~\ref{Lemm:TheoWAG-1}}.


\begin{figure}[t]
	\begin{center}
	\hspace{0mm}
    \includegraphics[width=150mm,clip]{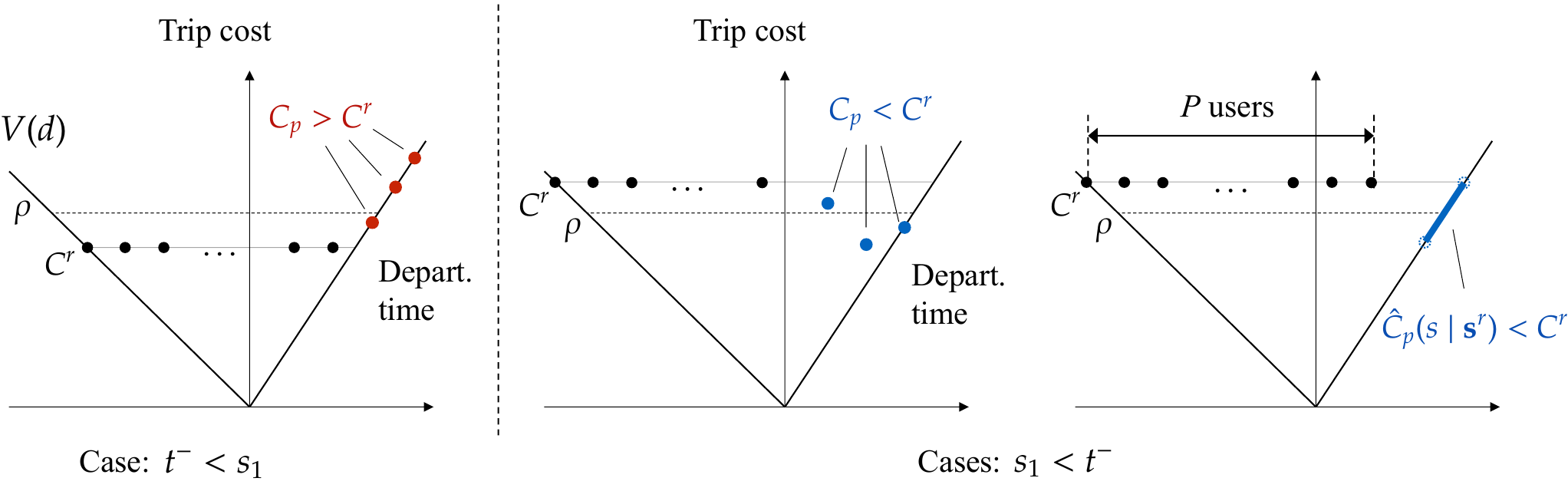}
	\end{center}
    \vspace{-4mm}
	\caption{Schematics of eventual time profiles when $C^{r}\neq \rho$.}
    \vspace{-2mm}
    \label{Fig:Reference}
\end{figure}

If the equilibrium cannot be achieved, Step 4 searches for the appropriate departure time of the first user, while narrowing the upper and lower bounds of the departure time range.
To adjust the departure time $s_{1}$ without a priori knowledge of $t^{-}$, this dynamics estimates whether $s_{1}$ is earlier or later than $t^{-}$ based on the eventual time profile where no more users can be fixed to the corresponding equilibrium departure time (see Figure~\ref{Fig:Reference}).
If $s_{1}$ is later than $t^{-}$, the shape of the schedule delay cost function implies that some users will experience trip costs higher than $C^{r}$ in the eventual time profile: the difference between the destination arrival times where the schedule delay costs are the same as $C^{r}$ is smaller than $m(P-1)/\mu$, which is the length of the rush-hour created by $P$ users.
Meanwhile, if $s_{1}$ is earlier than $t^{-}$, some users will experience trip costs lower than $C^{r}$, or all users will experience the same trip cost but the last user experiences a positive queueing delay cost.
Hence, it is expected to identify the direction of change to move $s_{1}$ close to $t^{-}$ from the eventual time profile.
After the estimation, the dynamics updates the upper or lower bounds of the departure time range and releases the fixed users to conduct better responses.
Users are fixed again after the departure time of the first user has been adjusted by better responses.

The following proposition ensures that one of the above cases indeed applies to the eventual time profile after better responses: 
\begin{prop}\label{Prop:AsymptoticFirst}
    Consider a time profile $\mathbf{s}$ where $s_{1} \neq t^{-}$.
    Suppose that the time profile satisfies the following conditions:
    (i) the first to $n$th users experience the same trip cost $C^{r}$ while arriving at the destination at the interval of $m/\mu$ and (ii) there does not exist a better response path to a time profile where $n+1$st user experiences $C^{r}$.
    
    Then, if $t^{-}<s_{1}$, there must exist users whose trip costs are different from $C^{r}$ and the trip costs of all such users are higher than $C^{r}$.
    If $s_{1}<t^{-}$, the time profile $\mathbf{s}$ satisfies one of the following situations: (a) all users experience the same trip cost $C^{r}$, but the last user experience a positive queueing delay cost; (b) there must exist users whose trip costs are different from $C^{r}$ and the trip costs of all such users are lower than $C^{r}$.

\end{prop}
\begin{prf}
See \textbf{\ref{Sec:App-Prop:AsymptoticFirst}}.\qed
\end{prf}

\noindent This proposition states that better responses with the fixation lead the time profile to a well-ordered traffic state although users' trip costs might be either higher or lower than $C^{r}$ during the transition process.\footnote{Figure~\ref{Fig:NonEquiSnapCost_500} in the numerical experiments shows that some users experience trip costs lower than $C^{r}$ even when the relationship $t^{-} < s_{1}$ holds.}
This also implies that the time profile does not get stuck into a non-equilibrium one where the departure time of the first user is not in equilibrium:
since users can conduct better responses in the eventual time profile as shown in Figure~\ref{Fig:Reference}, the departure time of the first user is adjusted by users' natural better responses once the fixation is released.

Combining \textbf{Theorem~\ref{Theo:WAG}} and \textbf{Proposition~\ref{Prop:AsymptoticFirst}}, we can establish the following global convergence of the better response dynamics to the epsilon-Nash equilibrium:
\begin{prop}\label{Prop:GlobalConvergence}
In a DTC game, the stationary point of the proposed better response dynamics is unique and corresponds to the epsilon-Nash equilibrium $\mathbf{s}^{*}$.
Moreover, a time profile almost surely converges to $\mathbf{r}^{*}$ from any time profile.
\end{prop}
\begin{prf}
See \textbf{\ref{Sec:App-Prop:GlobalConvergence}}.\qed



\end{prf}
\color{black}

\textcolor{black}{
This proposition states that we can construct the global convergent dynamics by incorporating the clarified convergence mechanisms.
Under the dynamics, the asymptotic adjustment process leads a traffic state to the exact starting point of the ordered path by changing the earliest departing user to the earliest equilibrium departure time.
Then, in the fixation process, users conduct better responses to the equilibrium departure times in the appropriate order.
As a result, a traffic state can converge to the equilibrium along the ordered path without knowing the path in advance.
}

\textcolor{black}{
A few remarks are in order.
First, }the global convergence implies the global stability because a time profile returns to the epsilon-Nash equilibrium $\mathbf{s}^{*}$ even if the time profile deviates from that equilibrium due to a perturbation.
In particular, the proposed evolutionary dynamics ensures that the departure time of the first user approaches the equilibrium departure time $t^{-}$;
more specifically, the difference in the upper and lower bounds of the departure time of the first user approaches zero while the equilibrium departure time is included between the bounds.

\textcolor{black}{
Second, the stability properties hold even when users are not bounded rational, i.e. the value of $\epsilon$ is extremely small.
As shown in Section 3, the epsilon ensuring the existence of the epsilon-Nash equilibrium approaches zero as the user size $m$ approaches zero.
Meanwhile, the stability results demonstrated in this section hold for an arbitrary $m$.
Therefore, even when $\epsilon$ is extremely small, we can still establish the stability results in a situation with a correspondingly small user size $m$.
This implies that the global convergence and its mechanisms do not vary even in a situation where the epsilon-Nash equilibrium~\eqref{Eq:DefiEpsilonNash} becomes almost the same as the pure Nash equilibrium~\eqref{Eq:DefiPureNash} by approaching the atomic model towards a fluid model.}

\textcolor{black}{
Third,} in Step 4 of the proposed dynamics, the global information (such as the time range) is used to adjust the departure time of the first user for the rigorous proof and numerical simulation of the global convergence of the dynamics to be shown later.
This is directly interpreted as being managed under a centralised control system; however, the adjustment might be also realised through decentralised better responses without using such global information.
This is because better responses in eventual time profiles will basically bring the departure time closer to the equilibrium one, as can be seen from Figure~\ref{Fig:Reference}.
Therefore, if the departure time of the first user gradually changes, the time profile is expected to converge to the equilibrium based on only the local information.
The development of more behaviourally natural dynamics based on the discussion is a subject for future research.

\color{black}

\subsection{Discussion about the stability and instability of other evolutionary dynamics}

We examine existing stable and unstable dynamics studied in the literature and discuss the reason for the different stability results based on the convergence mechanisms found in this study.
For stable dynamics, we examine the local dynamics of \cite{Jin2021-xe}.
This dynamics has the following two characteristics in ensuring stability: 
(i) limited available time periods: users can only change their departure times within the equilibrium rush hour $(t^{-},t^{+})$ and positive flows always occur within the time range; (ii) localness: users can only change their departure times to the adjacent departure times, e.g. from the time $s$ to $s+\Delta s$ or $s-\Delta s$.
Meanwhile, we consider other evolutionary dynamics, such as Smith dynamics, projection dynamics and replicator dynamics, for unstable dynamics.
They fundamentally share the following common features: 
the traffic flow at a departure time will shift to times that can improve the trip cost without any constraints on changeable times.

Referring to the convergence mechanisms, one of the important differences that characterise the stability between the two types of dynamics would be whether traffic flows are constrained to the range of the equilibrium rush hour.
This difference corresponds to whether the asymptotic adjustment mechanism of the earliest departing user is satisfied.
Under the stable dynamics, the characteristic of limited available time periods ensures that the departure time of the earliest traffic flow is fixed to the equilibrium departure time $t^{-}$.
This means that a traffic state belongs to a special initial time profile for convergence to the equilibrium.
Meanwhile, we easily see that each of the unstable dynamics lacks the asymptotic adjustment mechanism.
Namely, the range of available departure times is not limited to the equilibrium rush hour $[t^{-},t^{+}]$, and there is no force attracting the first user to the equilibrium departure time $t^{-}$.
The dynamics thus becomes unstable in the sense that it cannot trace the ordered path.

In addition, our convergence mechanisms suggest that the localness is also important to realise the stability.
This is because the localness works similarly to the fixation mechanism.
Suppose that traffic flows departing in a time range $[s^{-},s^{+}]$ are locally equilibrated while experiencing identical trip costs.
Then, the equilibrated traffic flows in $(s^{-},s^{+})$ have no incentive to change their departure times to the adjacent ones.
Furthermore, the equilibrated traffic flows are not perturbed by the departure time changes of other traffic flows departing far from that time range: such traffic flows cannot change their departure times to those in $(s^{-},s^{+})$ owing to the localness.
Therefore, we consider that the localness contributes to the local stability and its accumulation could lead to overall stabilisation, although this does not necessarily mean that the traffic flows become equilibrated starting from the earliest departure times, as in our study.

\color{black}

\section{Numerical experiments}
This section presents numerical experiments to demonstrate the convergence of the proposed better response dynamics to the epsilon-Nash equilibrium.
We first show the convergence from a special initial time profile for the convergence, in which the departure time of the first user is the same as $t^{-}$.
We next show that the dynamics can converge to the epsilon-Nash equilibrium even from a time profile other than the special one by asymptotic adjustment of the departure time of the first user.

\subsection{Settings}
We consider a single bottleneck with the capacity $\mu = 1.0$.
The number of users $P$ is $101$, and the user size $m$ is $1.0$.
The coefficients for the early and late arrival penalties are set as $\beta = 0.5$ and $\gamma = 2.0$, respectively.
Note that the desired arrival time and the free-flow travel time are $0$, as with the theoretical analysis.
Under these settings, the value of $\epsilon$, the equilibrium cost $\rho$ and the departure time of the first user in the epsilon-Nash equilibrium are $3.0$, $40.0$ and $-80.0$, respectively.
We also set $\Delta s$ as $0.01$.

Users change their departure times within the range of $[-100, 100]$, which includes the rush hour in the epsilon-Nash equilibrium, according to the proposed better response dynamics.
A selected user, who can change his/her departure time, searches for and selects a better response strategy randomly for the sake of simplicity.
Specifically, the selected user calculates the forecasted costs for a maximum of 100 random departure times that are later than the equilibrated users.
Then, if the user finds a departure time where the forecasted cost is lower than the current trip cost, the user changes to that departure time;
otherwise, the selected user does not change the departure time.
We assume that the selected user first checks the equilibrium (or reference) departure time, where the user can experience the same trip cost as that of the first user.
This is to accelerate the convergence to the equilibrium since the main aim of the numerical experiments is to confirm the convergence properties shown so far.


\begin{figure}
	\begin{minipage}[t]{0.33\textwidth}
		\centering
		\includegraphics[clip, width=1.0\columnwidth]{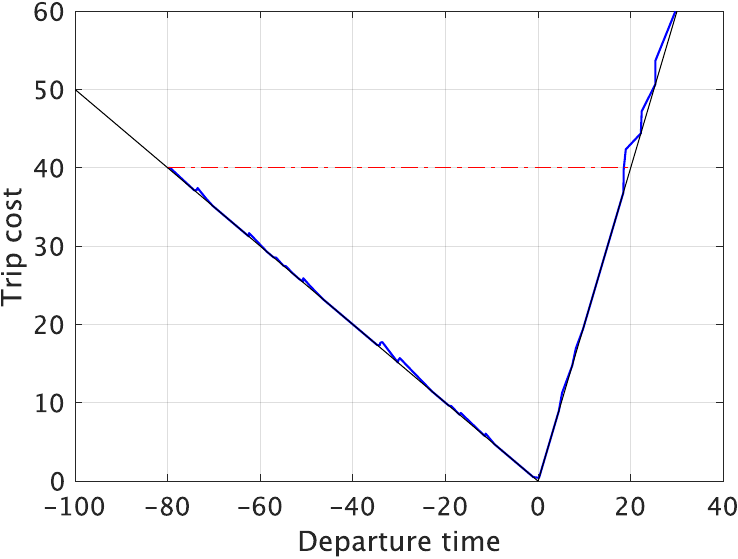}
		\subcaption{$\tau=1$ (Initial time profile)}\label{Fig:EquiSnapCost_1}
	\end{minipage}
	\begin{minipage}[t]{0.33\textwidth}
		\centering
		\includegraphics[clip, width=1.0\columnwidth]{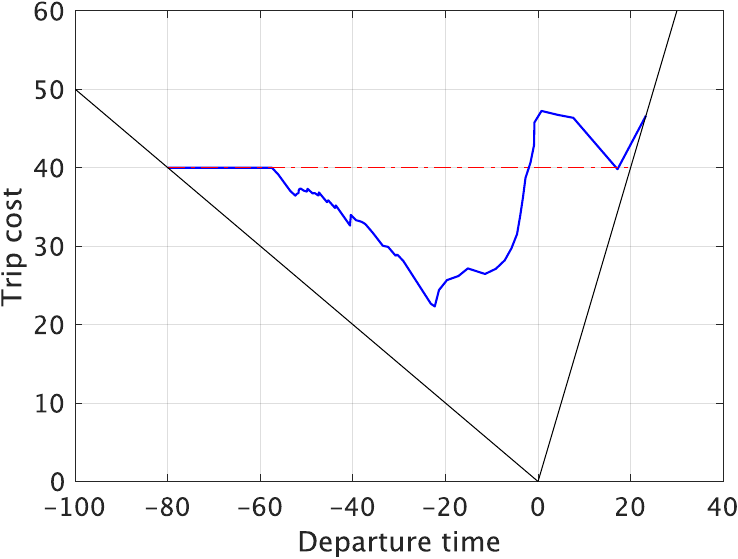}
		\subcaption{$\tau=500$}\label{Fig:EquiSnapCost_500}
	\end{minipage}
    \begin{minipage}[t]{0.33\textwidth}
		\centering
		\includegraphics[clip, width=1.0\columnwidth]{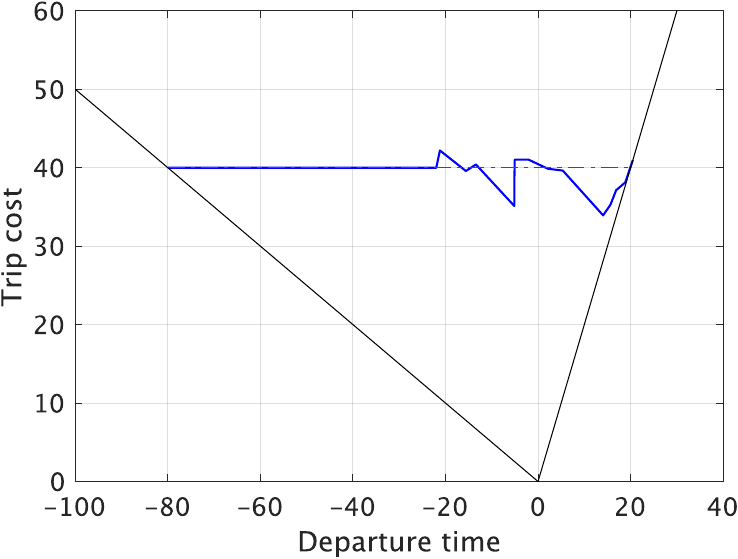}
		\subcaption{$\tau=900$}\label{Fig:EquiSnapCost_900}
	\end{minipage}
	\vspace{-1mm} 
	\caption{Snapshots of the trip costs of the users}
	\label{Fig:EquiSnapCost}
	\vspace{-2mm}
\end{figure}

\subsection{Results}

\subsubsection{Convergence from a special initial time profile}
\begin{figure}
    \begin{minipage}[t]{0.33\textwidth}
		\centering
		\includegraphics[clip, width=1.0\columnwidth]{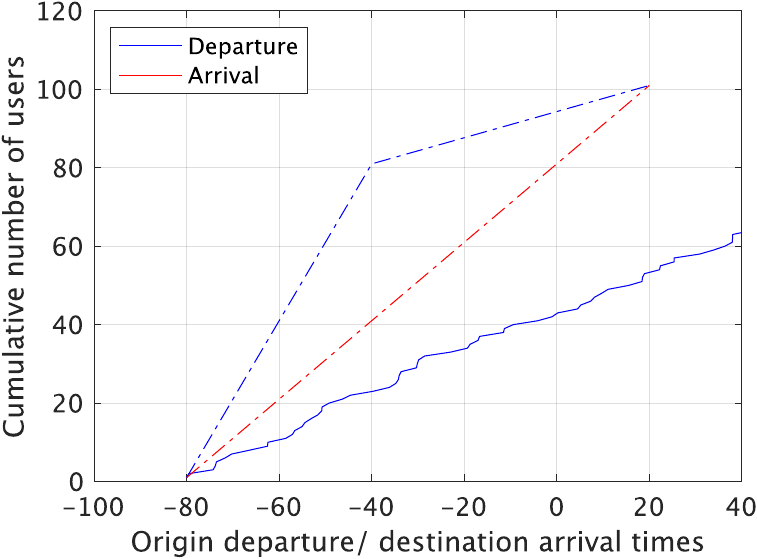}
		\subcaption{$\tau=1$}\label{Fig:EquiSnapFlow_1}
	\end{minipage}
	\begin{minipage}[t]{0.33\textwidth}
		\centering
		\includegraphics[clip, width=1.0\columnwidth]{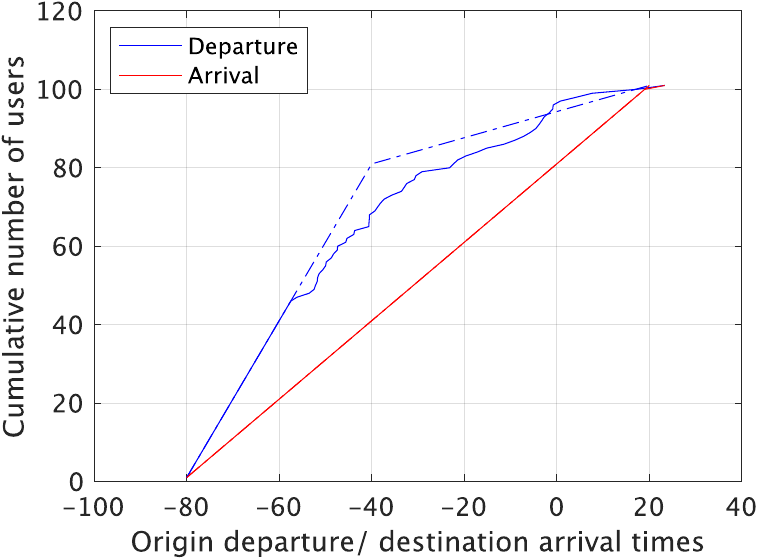}
		\subcaption{$\tau=500$}\label{Fig:EquiSnapFlow_500}
	\end{minipage}
    \begin{minipage}[t]{0.33\textwidth}
		\centering
		\includegraphics[clip, width=1.0\columnwidth]{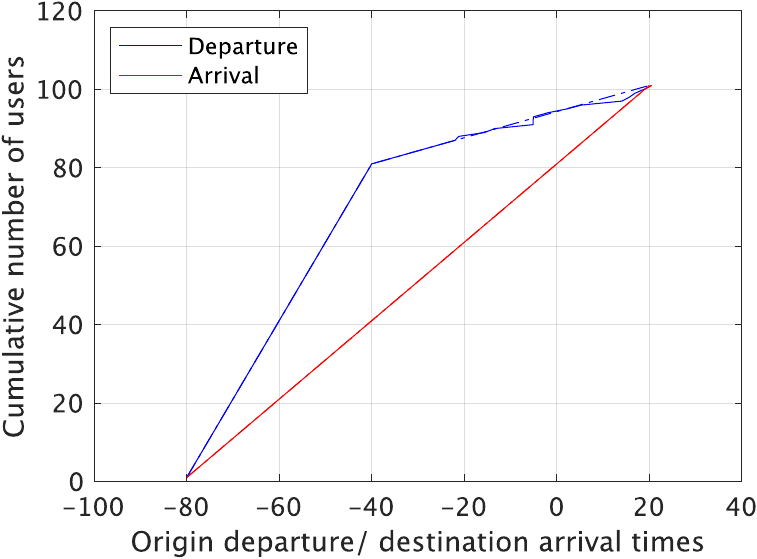}
		\subcaption{$\tau=900$}\label{Fig:EquiSnapFlow_900}
	\end{minipage}
	\vspace{-1mm} 
	\caption{Snapshots of the cumulative departure and arrival curves}
	\label{Fig:EquiSnapFlow}
	\vspace{-2mm}
\end{figure}

We first observe the convergence from an initial time profile where the first user departs at $t^{-}$.
The initial departure times of the other users are drawn from the uniform distribution between $t^{-} = -80$ and $100$.
Figure~\ref{Fig:EquiSnapCost} shows snapshots of the trip costs of the users at $\tau=1,500$ and $900$.
Figure~\ref{Fig:EquiSnapFlow} also shows the cumulative departure and arrival curves.
\textcolor{black}{For reference, the equilibrium trip cost and departure/arrival curves are plotted as dash-dotted lines.}
We confirm from these figures that users' trip costs and departure times are sequentially fixed to those in the epsilon-Nash equilibrium.
We also see that equilibrated users do not change their departure times, and they are not perturbed by the overtaking of non-equilibrated users.

\begin{figure}[t]
	\begin{center}
	\hspace{0mm}
    \includegraphics[width=160mm,clip]{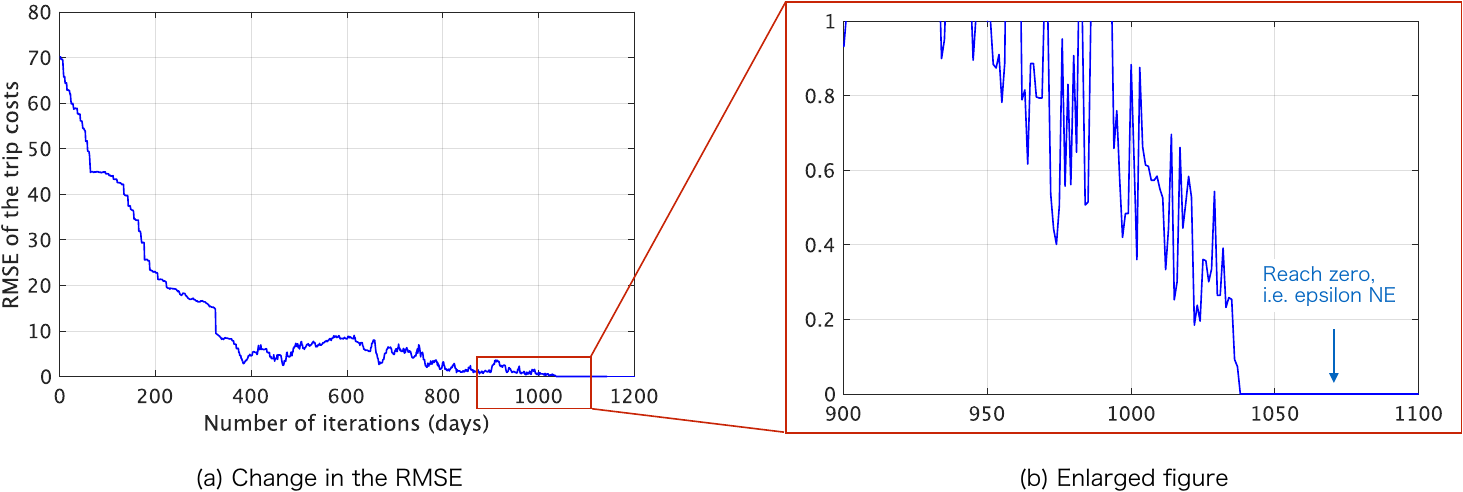}
	\end{center}
    \vspace{-4mm}
	\caption{Convergence process of the better response dynamics}
    \vspace{-0mm}
    \label{Fig:EquiCostChange}
\end{figure}

Figure~\ref{Fig:EquiCostChange} shows the root mean squared error (RMSE) of the trip costs for the equilibrium cost $\rho$ on each day $\tau$, which is defined as follows:
\begin{align}
    z^{\tau} = \left[\cfrac{1}{P}\sum_{p\in\mathcal{P}}(C_{p}(\mathbf{s}^{\tau}) - \rho)^{2} \right]^{1/2}.
\end{align}
\textcolor{black}{The figure indicates that the RMSE reaches zero, which means that the time profile converges to the epsilon-Nash equilibrium.}
These results numerically confirm the existence of the better response path leading to the equilibrium, as mentioned in \textbf{Section~\ref{Sec:WAG}}.

\begin{figure}
	\begin{minipage}[t]{0.33\textwidth}
		\centering
		\includegraphics[clip, width=1.0\columnwidth]{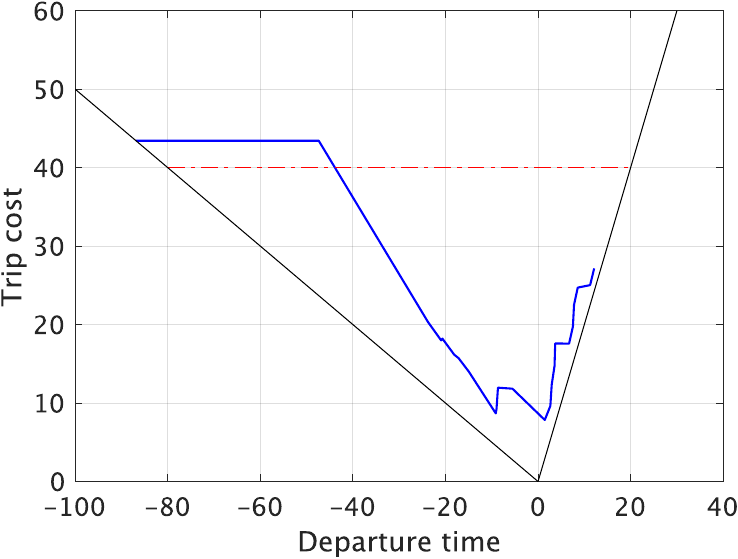}
		\subcaption{$\tau=98000$}\label{Fig:NonEquiSnapCost_1}
	\end{minipage}
	\begin{minipage}[t]{0.33\textwidth}
		\centering
		\includegraphics[clip, width=1.0\columnwidth]{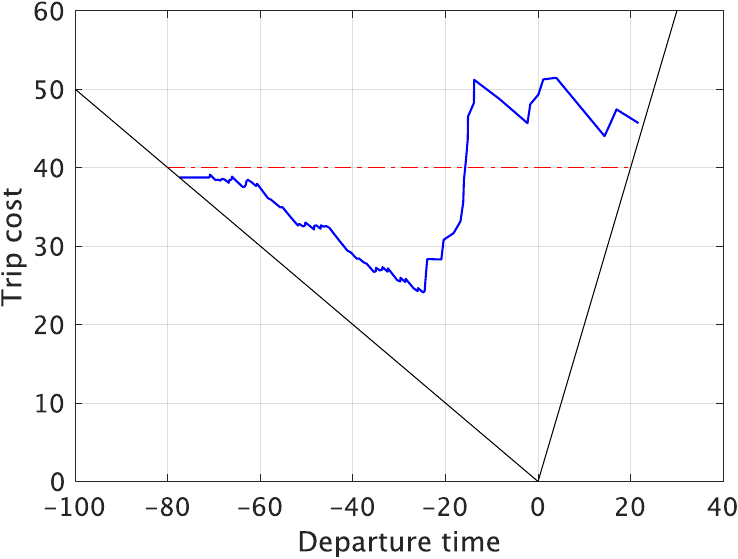}
		\subcaption{$\tau=160000$}\label{Fig:NonEquiSnapCost_500}
	\end{minipage}
    \begin{minipage}[t]{0.33\textwidth}
		\centering
		\includegraphics[clip, width=1.0\columnwidth]{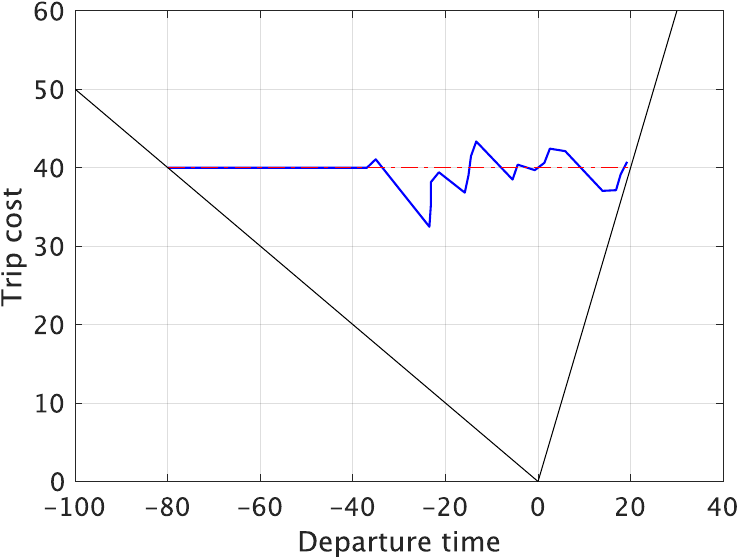}
		\subcaption{$\tau=248800$}\label{Fig:NonEquiSnapCost_900}
	\end{minipage}
	\vspace{-1mm} 
	\caption{Snapshots of trip costs of the users: general initial time profile case}
	\label{Fig:NonEquiSnapCost}
	\vspace{-2mm}
\end{figure}

\subsubsection{Convergence from a general initial time profile}
We next examine the behaviour of the better response dynamics from a time profile where the first user does not depart at $t^{-}$.
In this case, for a fixed departure time of the first user, a time profile may be led to a situation where all of the users cannot conduct better responses such that they experience the trip cost $C^{r}$, which is the same as that of the first user.
To avoid getting stuck in such a time profile, we check the time profile if the number of users who experience $C^{r}$ does not increase even after $10,000$ better responses.
We then update the upper and lower bounds for the departure time of the first user, as mentioned in the previous section.

\begin{figure}
    \begin{minipage}[t]{0.33\textwidth}
		\centering
		\includegraphics[clip, width=1.0\columnwidth]{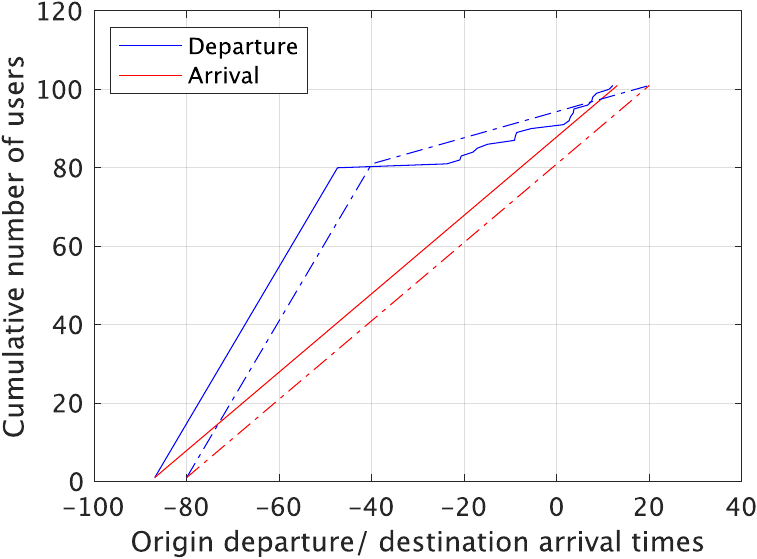}
		\subcaption{$\tau=98000$}\label{Fig:NonEquiSnapFlow_30000}
	\end{minipage}
	\begin{minipage}[t]{0.33\textwidth}
		\centering
		\includegraphics[clip, width=1.0\columnwidth]{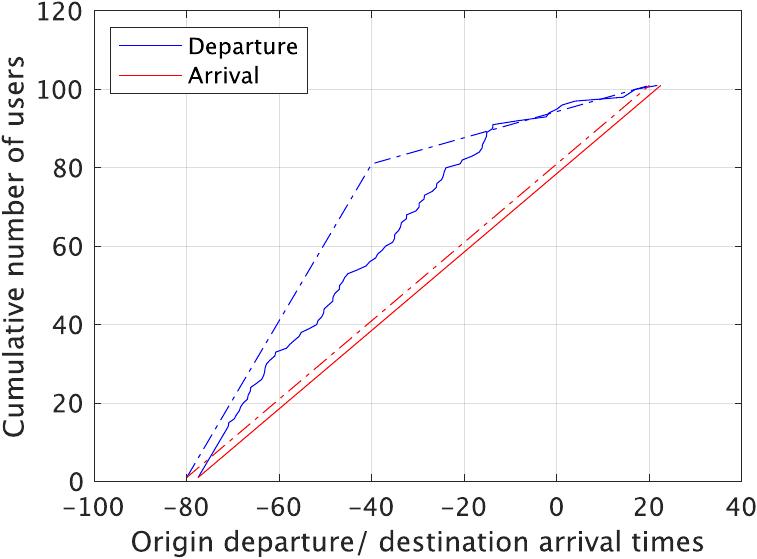}
		\subcaption{$\tau=160000$}\label{Fig:NonEquiSnapFlow_60000}
	\end{minipage}
    \begin{minipage}[t]{0.33\textwidth}
		\centering
		\includegraphics[clip, width=1.0\columnwidth]{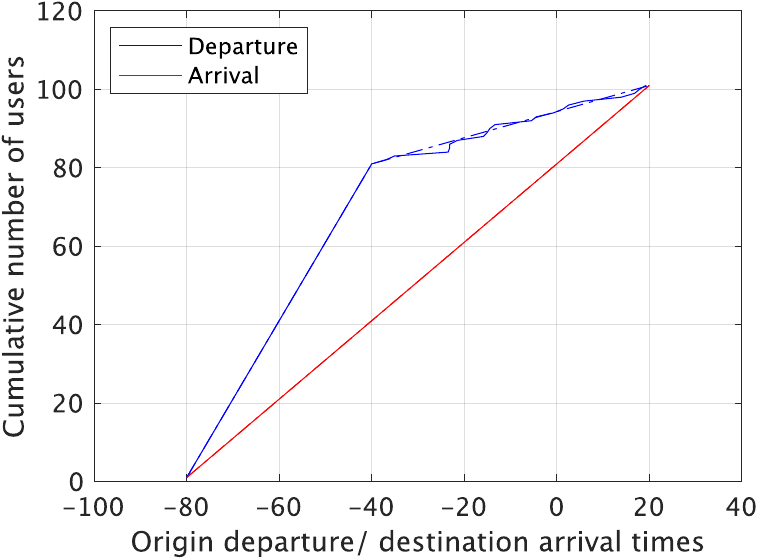}
		\subcaption{$\tau=248800$}\label{Fig:NonEquiSnapFlow_7000}
	\end{minipage}
	\vspace{-1mm} 
	\caption{Snapshots of cumulative departure/arrival curves: general initial time profile case}
	\label{Fig:NonEquiSnapFlow}
	\vspace{-2mm}
\end{figure}

Figure~\ref{Fig:NonEquiSnapCost} shows snapshots of the trip costs of the users; Figure~\ref{Fig:NonEquiSnapFlow} shows snapshots of the cumulative curves.
\textcolor{black}{The dash-dotted lines represent the equilibrium trip cost and departure/arrival curves.}
These figures suggest the convergence of the departure time of the first user to the equilibrium one without knowing it in advance.
Specifically, when the departure time is earlier than $t^{-}$, the departure time changes so as to approach $t^{-}$.
When the departure time becomes later than $t^{-}$, the departure time changes to an earlier time.
Then, the departure time finally succeeds in changing to $t^{-}$ through better responses.

\begin{figure}[t]
	\begin{center}
	\hspace{0mm}
    \includegraphics[width=70mm,clip]{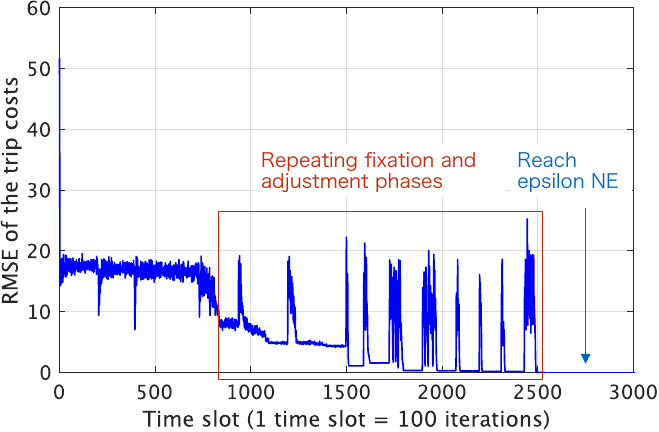}
	\end{center}
    \vspace{-4mm}
	\caption{Convergence process: general initial time profile case}
    \vspace{-2mm}
    \label{Fig:NonEquiCostChange}
\end{figure}

Figure~\ref{Fig:NonEquiCostChange} shows the RMSE of the time profile on each day $\tau$ to the equilibrium.
From this figure, we can confirm that the time profile converges to the epsilon-Nash equilibrium while oscillating due to the adjustment of the first user's departure time and fixation of the subsequent users.
In addition, we see that the RMSEs tend to decrease with changing the departure time of the first user (after jumps in RMSEs).
This implies that the time profile approaches the epsilon-Nash equilibrium on average as the departure time of the first user approaches the equilibrium one.
In summary, we conclude that the evolutionary dynamics can globally converge to the epsilon-Nash equilibrium without depending on the initial time profile.

\section{Conclusion}
This study conducted a stability analysis of equilibrium in DTC problems using an atomic approach.
First, we formulated a DTC game, which is a strategic game form of the DTC problems in which atomic users select their departure times.
We introduced the concept of the epsilon-Nash equilibrium as the equilibrium of this game.
This enables us to obtain rigorous results regarding the existence of pure-strategy equilibrium corresponding to the equilibrium of conventional fluid models.
We then proved that the DTC game is a weakly acyclic game in the sense that there exists a better response path converging to the epsilon-Nash equilibrium from an arbitrary initial traffic state.
Such a convergent better response path allowed us to clarify the mechanisms for global convergence to the epsilon-Nash equilibrium.
Based on the convergence mechanisms, we established the better response dynamics under which the equilibrium became globally stable.
We also discussed the underlying factors of the existing stable and unstable dynamics that yielded different stability results based on these mechanisms.
Numerical experiments confirmed that the proposed better response dynamics converged globally to the epsilon-Nash equilibrium.

While we presented the evolutionary dynamics that globally converged to equilibrium, this dynamics is only a naive formulation of sufficient conditions to guarantee the convergence.
Therefore, it is important to develop behaviourally natural evolutionary dynamics.
For this issue, we can consider not only myopic evolutionary dynamics but also dynamics incorporating natural individual behaviour (e.g. learning behaviour), which can be implemented using our atomic approach straightforwardly.
Moreover, it would be interesting to analyse the relationship between the macroscopic behaviour of traffic flow and the microscopic behaviour of each individual user.
This can contribute to the development of stable dynamics in fluid models, for example, through the deterministic approximation of the stochastic process (i.e. Markov chain induced from evolutionary dynamics) of individual DTCs~\citep[][]{Sandholm2010-ht,Iryo2016-pj}.

\textcolor{black}{
It would also be worthwhile to extend our analysis to more general frameworks of DTC problems with bottleneck models, such as with multiple bottlenecks, competing transport modes and user heterogeneity.
Because the FIFO principle and causality properties still hold in such extended models, the stability could be analysed based on the ordering property in a manner similar to this study.
Another important challenge is to incorporate hypercongestion phenomena of Network Fundamental Diagram (NFD) models~\citep[e.g.][]{Arnott2013-oq,Jin2020-mt}.
Regarding this topic, the atomic model seems to have a similar mathematical structure to the trip-based NFD models in the sense that both models track the trip (length) of individual users; 
indeed, the NFD models have recently been discussed in relation to an agent-based model~\citep[e.g.][]{Daganzo2015-zb,Mariotte2017-yc,Lamotte2018-ib}. 
Therefore, we expect that it may be possible to formulate a new DTC game by incorporating the NFD model, and conduct a stability analysis based on game theory, as shown in this study.
Exploring the stability in these extended situations from a game-theoretic perspective is an important topic for future research.
}

\section*{Declaration of competing interest}
None.

\section*{Acknowledgement}
We would like to thank the anonymous referees for their constructive comments, which helped to improve this paper.
This work was supported by JSPS KAKENHI Grant Numbers JP23K13418, JP23H01524 and JP20H00265.

\appendix

\section{Proof of the existence of the epsilon-Nash equilibrium}\label{App:Prf_Improvement}

We first establish the following lemma regarding the property of the schedule delay cost function:
\begin{lemm}\label{Lemm:App-RelationSchedule}
Consider a real number $t\in\mathbb{R}$ and a non-negative real number $\Delta t\in\mathbb{R}^{+}_{0}$.
Then, the following relationship holds true:
\begin{align}
    -\beta \Delta t \leq V(t+\Delta t) - V(t) \leq \gamma \Delta t.
\end{align}
\end{lemm}
\begin{prf}
Suppose first the case where $t+\Delta t\leq 0$ holds.
We then have
\begin{align}
    V(t+\Delta t) - V(t) = -\beta (t+\Delta t) + \beta t = -\beta \Delta t.
\end{align}
Suppose next the case where $t\geq 0$ holds.
We then have
\begin{align}
    V(t+\Delta t) - V(t) = \gamma (t+\Delta t) - \gamma t = \gamma \Delta t.
\end{align}
Finally, suppose the case where $0<t+\Delta t$ and $t<0$.
It follows that $-\Delta t < t < 0$.
We then have
\begin{align}
    &V(t+\Delta t) - V(t) = (\beta + \gamma)t + \gamma \Delta t, \quad \Rightarrow \quad -\beta \Delta t < V(t+\Delta t) - V(t) < \gamma \Delta t.
\end{align}
These results prove the lemma.\qed
\end{prf}

Utilising this lemma, we calculate the maximum improvement in trip cost achievable by one user changing his/her departure time unilaterally in the equilibrium state $\mathbf{s}^{e}$.
Consider a situation where a user $p\in\mathcal{P}$, whose order of departure is $o_{p}$, changes the departure time to $s'_{p}$.
We denote by $o'_{p}$, $d'_{p}$ and $C'_{p}$ the order of departure, destination arrival time and trip cost of the user in the new time profile $\mathbf{s}'\equiv (s'_{p},\mathbf{s}_{-p}^{e})$, respectively.
For the sake of simplicity, we sometimes use the notation $d_{p}$ to represent the destination arrival time of the user in $\mathbf{s}^{e}$: i.e. $d_{p} = d_{o_{p}}$.
We then calculate the maximum improvement by dividing the analysis into the following three cases: (i) $o'_{p} < o_{p}$, (ii) $o_{p} = o'_{p}$ and (iii) $o_{p} < o'_{p}$.

\subsection{Case (i)}
In this case, user $p$ changes the departure time while moving up the departure order, which means that there exists a user behind user $p$ after the change.
This case further implies that we do not need to consider the case where $o'_{p} = 1$, i.e. user $p$ changes the departure time to depart earliest: in this situation, since $s'_{p}$ becomes earlier than $t^{-}$,  the schedule delay cost alone exceeds $\rho$, and thus the user cannot improve the trip cost.
It is thus sufficient for us to consider a situation where there exist users in front of and behind user $p$ after the change.

Let $p_{f}$ and $p_{b}$ denote the users in front of and behind user $p$ in $\mathbf{s}'$, respectively.
It follows that $s'\in(s_{p_{f}},s_{p_{b}})$.
In addition, $d_{o_{p_{b}}} - d_{o_{p_{f}}} = m/\mu$ holds in $\mathbf{s}^{e}$ and $d_{o_{p_{f}}}$ does not change by the departure time change of user $p$ who departs later than user $p_{f}$.
Therefore, user $p$ who overtakes the user $p_{b}$ from behind must travel in a congested situation in $\mathbf{s}'$, which means that 
\begin{align}
    d'_{p} = d_{o_{p_{f}}} + m/\mu.
\end{align}
We then obtain the trip cost of user $p$ as follows:
\begin{align}
C'_{p} 
    &= ( d'_{p}- s'_{p} ) + V(d'_{o'_{p}})\notag\\
    &= ( d_{o_{p_{f}}} -s_{p_{f}}) + V(d_{o_{p_{f}}}) + (s_{p_{f}} - s'_{p} )  + d'_{p} - d_{o_{p_{f}}} + V(d_{o_{p_{f}}}+m/\mu) -  V(d_{p_{f}})\notag\\
    &=\rho + (s_{p_{f}} - s'_{p} )  + m/\mu + V(d_{o_{p_{f}}}+m/\mu) -  V(d_{o_{p_{f}}}).
\end{align}
Moreover, using \textbf{Lemma~\ref{Lemm:App-RelationSchedule}} and the fact that $s'_{p}<s_{p_{b}} \leq s_{p_{f}}+m(1+\gamma)/\mu$, we have
\begin{align}
    C'_{p}\geq \rho + (s_{p_{f}} - s'_{p} ) +\cfrac{m(1-\beta)}{\mu}
    > \rho - \cfrac{m(\beta + \gamma)}{\mu}.
\end{align}
Then, the improvement in the trip cost of the user $p$ is described as follows:
\begin{align}
    C_{p} - C'_{p} < \cfrac{m(\beta + \gamma)}{\mu} <\cfrac{m(1+\gamma)}{\mu}.
\end{align}

\subsection{Case (ii)}
In this case, user $p$ changes the departure time while keeping the departure order.
We first describe the trip cost $C'_{p}$ as follows:
\begin{align}
    C'_{p} &= (d'_{p} - s'_{p}) + V(d'_{p}) \notag\\
    &= (d_{p} - s_{p}) + V(d_{p}) - (d_{p} - s_{p}) - V(d_{p})  + (d'_{p} - s'_{p}) + V(d'_{p}) \notag\\
    &= \rho + (d'_{p} - d_{p}) + (s_{p} - s'_{p}) + V(d'_{p}) - V(d_{p}).
\end{align}
We then divide the analysis into the following cases: (a) $1<o_{p}<P-1$, (b) $o_{p} = P$ and (c) $o_{p} = 1$.

In Case (a), $s' < s_{p}+m(1+\gamma)/\mu$.
In addition, since $d_{o_{p}} - d_{o_{p} - 1} = m/\mu$, the destination arrival time $d'_{p}$ cannot become earlier than $d_{p}$.
Thus, it can be described by using non-negative real number $\Delta d\geq 0$ as follows: 
\begin{align}
    d'_{p} = d_{p} + \Delta d.\label{Eq:App-Epsilon_Caseii-a-1}
\end{align}

Then, $C'_{p}$ is described as follows:
\begin{align}
    C'_{p} &\geq \rho + (s_{p} - s'_{p}) + \Delta d(1-\beta) \geq \rho + (s_{p} - s'_{p}).
\end{align}
We thus have
\begin{align}
    C_{p} - C'_{p} \leq s'_{p} - s_{p} < \cfrac{m(1+\gamma)}{\mu}.
\end{align}

In Case (b), the destination arrival time $d'_{p}$ can be described as Eq.~\eqref{Eq:App-Epsilon_Caseii-a-1} since the user also travels in a congested situation in $\mathbf{s}^{e}$.
In addition, since we consider the last user, the following relationships hold true:
\begin{align}
    s_{p} = d_{p} = t^{+},\quad V(d_{p}) = \gamma d_{p} = \rho.
\end{align}
We thus have
\begin{align}
    C'_{p}\geq V(d'_{p}) \geq V(d_{p}) = \rho,\quad \Rightarrow \quad C_{p} - C'_{p}\leq 0.
\end{align}

In Case (c), $s'_{p}\leq s_{p}+m(1+\gamma)/\mu$.
When $o_{p} = 1$, we do not need to consider a case where the destination arrival time becomes earlier: it is obvious that $d'_{p}<t^{-}$ and thus the schedule delay cost alone exceeds $\rho$.
Thus, the destination arrival time $d'_{p}$ can be described as Eq.~\eqref{Eq:App-Epsilon_Caseii-a-1}.
Since the analysis condition is the same as Case (a), we have 
\begin{align}
    C_{p} - C'_{p} \leq s'_{p} - s_{p} \leq \cfrac{m(1+\gamma)}{\mu}.
\end{align}

\subsection{Case (iii)}
In this case, user $p$ changes their departure time while delaying the order of departure, which means that there exists a user in front of user $p$ after the change.
This case further implies that we do not need to consider the case where $o'_{p} = P$, i.e. user $p$ changes the departure time so as to depart latest: in this situation, since $s'_{p}$ becomes later than $t^{+}$,  the schedule delay cost alone exceeds $\rho$ and thus the user cannot improve the trip cost.
It is thus sufficient for us to consider a situation where there exist users in front of and behind user $p$ after the change.

Let $p_{f}$ and $p_{b}$ denote the users in front of and behind user $p$ in $\mathbf{s}'$, respectively.
We also denote by $d'_{p_{f}}$ the destination arrival time of user $p_{f}$ in $\mathbf{s}'$.
This arrival time can be at most $m/\mu$ earlier than $d_{p_{f}}$, which is the arrival time of the user in $\mathbf{s}$ because there is one less user travelling in front of the user in $\mathbf{s}'$.
Thus, the arrival time can be described as follows:
\begin{align}
    d'_{p_{f}} + \Delta n = d_{p_{f}},
\end{align}
where $\Delta n$ is a non-negative number satisfying $0\leq \Delta n \leq m/\mu$.
By using this variable, the destination arrival time $d'_{p}$ of user $p$ after the change can be described as follows:
\begin{align}
    d'_{p} = d'_{p_{f}} + \cfrac{m}{\mu} + \Delta k,
\end{align}
where $\Delta k$ is a non-negative number, i.e. $\Delta k\geq 0$.

Then, the trip cost $C'_{p}$ can be described as follows:
\begin{align}
C'_{p} 
&= (d'_{p} - s'_{p}) + V(d'_{p}) \notag \\
&= (d_{p_{f}} - s_{p_{f}}) + V(d_{p_{f}}) - (d_{p_{f}} - s_{p_{f}}) - V(d_{p_{f}})  + (d'_{p} - s'_{p}) + V(d'_{p}) \notag \\
&= \rho  - (\Delta n - \Delta k - m/\mu) +  (s_{p_{f}} - s'_{p}) + V(d'_{p}) - V( d'_{p} + \Delta n- \Delta k - m/\mu  ).
\end{align}

\noindent Here, we have the following relationship:
\begin{align}
    \Delta n - \Delta k - m/\mu\leq -\Delta k \leq 0.
\end{align}
We thus have
\begin{align}
    V(d'_{p}) - V( d'_{p} + \Delta n - \Delta k - m/\mu  ) \geq \beta ( \Delta n -\Delta k - m/\mu  ).
\end{align}
Substituting this, we have
\begin{align}
C'_{p} 
&\geq  \rho +  (s_{p_{f}} - s'_{p}) - (1-\beta)( \Delta n - m/\mu - \Delta k ) \geq \rho +  (s_{p_{f}} - s'_{p}).
\end{align}
Thus, the improvement is represented as follows:
\begin{align}
    C_{p} - C'_{p} \leq s'_{p} - s_{p_{f}} < \cfrac{m(1+\gamma)}{\mu}.
\end{align}

In summary, we see that an improvement in the trip cost of a user does not exceed the value of $m(1+\gamma)/\mu$.
Therefore, by setting $\epsilon\geq m(1+\gamma)/\mu$, any user cannot improve his/her trip cost by more than $\epsilon$ by unilaterally changing his/her departure time in time profile $\mathbf{s}^{e}$.
This means that the time profile becomes $\epsilon$-Nash equilibrium, and this proves \textbf{Theorem~\ref{Theo:ExistenceNash}}.\qed

\section{Setting of $\Delta s$}\label{Sec:App-SettingDelta_t}
\begin{lemm}\label{Lemm:App-Delta_t-1}
Consider a time profile $\mathbf{s}^{e}$, where all users experience the same trip cost $\rho$.
Then, the discretisation value $\Delta s$ satisfies the following conditions such that the users can select the equilibrium departure time precisely:
\begin{align}
    \cfrac{\rho}{\beta} = k^{-}\cdot \Delta s,\quad
    \cfrac{m(1-\beta)}{\mu} = k^{\beta}\cdot \Delta s,\quad
    \cfrac{m(1+\gamma)}{\mu} = k^{\gamma}\cdot \Delta s,\quad 
    (\beta+\gamma)k^{-}\Delta s = k^{+}\cdot \Delta s,
\end{align}
where $k^{-}$, $k^{\beta}$, $k^{\gamma}$ and $k^{+}$ are positive integers, i.e. $\Delta s$ is set such that these values are multiples of $\Delta s$.
\end{lemm}
\begin{prf}
    First, $\rho/\beta$ is equivalent to $t^{-}$ and thus this must be a multiple of $\Delta s$.
    Second, $m(1-\beta)/\mu$ and $m(1+\gamma)/\mu$ are the intervals between departure times of early and late arrival users, respectively; thus, these must be multiples.
    
    Finally, the fourth condition is needed to ensure that the interval between the departure time of the latest user $p$ among early arrival users (i.e. $o_{p}=o^{cr}$) and the earliest user $p'$ among late arrival users (i.e. $o_{p'}=o^{cr}+1$) is a multiple of $\Delta s$.
    The interval is calculated as follows:
    \begin{align}
        s_{p'}^{e} - s_{p}^{e} = \left(\cfrac{m(1+\gamma)}{\mu} - \cfrac{m(1-\beta)}{\mu}\right)\left\lfloor \cfrac{\gamma(P-1)}{\beta + \gamma}\right\rfloor
        +\cfrac{m(1+\gamma)}{\mu}
        -\cfrac{(\beta + \gamma)\rho}{\beta}
    \end{align}
    The first and second terms are obviously multiples of $\Delta s$.
    Then, the third term must be a multiple of $\Delta s$ and this holds when $(\beta + \gamma)k^{-}\Delta s$ is a multiple of $\Delta s$.
    Note that such an integer number $k^{+}$ always exists when $\beta+\gamma$ is a rational number and this is ensured since $\beta$ and $\gamma$ are rational numbers.\qed
\end{prf}

\begin{lemm}
Suppose that $\Delta s$ is set such that $m/\mu$ is a multiple of $\Delta s$, i.e. $m/\mu$ can be represented as $k^{\mu}\cdot \Delta s$.
Then, the destination arrival time of any user becomes a multiple of $\Delta s$.
\end{lemm}
\begin{prf}
Consider first a user $p\in\mathcal{P}$ who travels in a free-flow situation.
Since the destination arrival time becomes the departure time which is a multiple of $\Delta s$, the statement is true.

Consider next a user $p\in\mathcal{P}$ who travels in a congested situation.
Then, the relationship $d_{o_{p}} - d_{o_{p}-1} =m/\mu$.
This means that $d_{o_{p}}$ is a multiple of $\Delta s$ if $d_{o_{p}-1}$ is also a multiple of $\Delta s$.
Then, by recursively applying this relationship, we find a user $p'$ who travels in a free-flow situation.
Therefore, $d_{o_{p}}$ is also a multiple of $\Delta s$.\qed
\end{prf}

\begin{lemm}\label{Lemm:App-Delta_t-2}
Suppose that $m/\mu$ is a multiple of $\Delta s$.
Consider two variables representing departure or destination arrival times $t$ and $t'$.
Then, the following relationship holds true.
\begin{align}
    t' > t\quad \Rightarrow \quad  t' \geq t+\Delta s.
\end{align}
\end{lemm}
\begin{prf}
Since these variables are multiples of $\Delta s$, the difference must be a multiple of $\Delta s$.
This proves the relationship.\qed
\end{prf}

\noindent In the rest of the appendices, which include \textbf{\ref{App:Prf_WAG}} and \textbf{\ref{Sec:App-Prop:AsymptoticFirst}}, we assume that $\Delta s$ is sufficiently small such that $m/\mu$ is a multiple of the value and satisfies the conditions in \textbf{Lemma~\ref{Lemm:App-Delta_t-1}}.

\section{Proof of the weakly acyclicity of DTC games}\label{App:Prf_WAG}
This section proves the two lemmas in the proof of \textbf{Theorem~\ref{Theo:WAG}}.
Hereinafter, we accord the index of a user with the order of departure.
For example, a notation $s_{k}$ denotes the departure time of the $k$th user departing from the origin in time profile $\mathbf{s}$.

Prior to the proof, we formally describe the numerical expression of the forecasted cost.
Consider a departure time $s'$ between $s_{k}$ and $s_{k+1}$ $(1\leq k\leq P-1)$ when the current time profile is $\mathbf{s}$.
Then, the forecasted trip cost, $\hat{C}_{p}(s' \mid \mathbf{s})$, is defined as follows (we hereinafter omit the subscript $p$ since it is irrelevant who changes the departure time in the definition under the assumption of the homogeneity of users):
\begin{itemize}
    \item if $d_{k+1} - d_{k} = m/\mu$:
    \begin{align}
        \hat{C}(s' \mid \mathbf{s}) = 
        C_{k} + \cfrac{C_{k+1} - C_{k}}{s_{k+1} - s_{k} }(s'-s_{k}).\label{Eq:DefiBetterCongested}
    \end{align}

    \item otherwise:
    \begin{align}
        \hat{C}(s' \mid \mathbf{s}) = 
        \begin{cases}
            C_{k} + 
            \cfrac{V(d_{k}) - C_{k}}{d_{k} - s_{k}}(s'-s_{k}),
            &\quad \text{if}\quad s'\in(s_{k},d_{k} ],\\
            V(s')
            &\quad \text{if}\quad s'\in(d_{k},s_{k+1}).\\
        \end{cases}\label{Eq:DefiBetterFree}
    \end{align}
\end{itemize}
Note that if $s' < s_{1}$ or $s' \geq d_{P}$, then $\hat{C}(s' \mid \mathbf{s}) = V(s')$.
Below, we omit $\mathbf{s}$ from $\hat{C}(s' \mid \mathbf{s})$ when the specific current time profile $\mathbf{s}$ is clear.

\subsection{Preliminaries}
We here establish several lemmas.
We first prove the following lemma regarding the relationship between the trip cost of the last user departing from the origin:
\begin{lemm}\label{Lemm:App-FirstLast}
Consider an arbitrary time profile $\mathbf{s}$.
Suppose that $t^{-}\leq s_{1}$ holds.
Then, 
\begin{align}
    C_{P}(\mathbf{s})\geq \rho.\label{Eq:App-FirstLastCost}
\end{align}
If $t^{-} < s_{1}$, the relationship holds strictly, i.e. $C_{P}>\rho$.
\end{lemm}
\begin{prf}
    Eq.~\eqref{Eq:QueueCost} suggests that the interval between the destination arrival times of two consecutive users is equal to or more than $m/\mu$.
    Since the first user is followed by $(P-1)$ users, $d_{P}\geq s_{1}+m(P-1)/\mu$ holds true.
    We thus have
    \begin{align}
        &t^{+} = t^{-}+\cfrac{m(P-1)}{\mu} \leq  s_{1}+\cfrac{m(P-1)}{\mu} \leq d_{P},\quad \Rightarrow \quad 
        C_{P}(\mathbf{s})\geq \gamma d_{P} \geq  \gamma t^{+} = \rho.
    \end{align}
    It is obvious that the relationship holds strictly if $t^{-} < s_{1}$.\qed
\end{prf}

We next establish the following lemma that states the existence of a departure time when a user can travel in a free-flow situation in some specific time profiles:
\begin{lemm}\label{Lemm:App-FreeFlowRange}
Consider a time profile $\mathbf{s}$.
Suppose that the relationship $d_{P}-s_{1}>m(P-1)/\mu$ is satisfied.
Then, there exists an appropriate integer number $n$ $(0\leq n \leq P)$ such that there exists a departure time $s_{F}$ satisfying the following relationship:
	\begin{align}
        d_{n} + \cfrac{m}{\mu}  \leq s_{F} < d_{n+1},
    \end{align}
    where $d_{0} = -\infty$ and $d_{P+1} = \infty$.    
\end{lemm}
\begin{prf}
Since $t^{+}-t^{-}=m(P-1)/\mu$, all $P$ users must arrive at the destination at the interval of $m/\mu$ in $[t^{-},t^{+}]$ if such a departure time does not exist, i.e. $d_{1}=t^{-}$ and $d_{P}=t^{+}$.
However, this contradicts the supposition.
Thus, the lemma is true.
\qed
\end{prf}

\noindent We have the following corollary from the lemma:
\begin{coro}\label{Coro:App-FreeFlowRange}
Consider a time profile $\mathbf{s}$ satisfying the relationship $d_{P}-s_{1}>m(P-1)/\mu$.
Then, there exists a departure time where the forecasted cost is $V(s_{F})\leq \rho$.
In addition, when a user departing earlier than $s_{F}$ changes the departure times to $s_{F}$, the destination arrival times of the users arriving at the destination later than $s_{F}$ before the change do not become earlier by the change.
\end{coro}

\color{black}
We further establish the following two lemmas that become the keys to proving the existence of an ordered path:
\begin{lemm}\label{Lemm:SnR<Sn}
Consider a time profile $\mathbf{s}$ where the first and $n$th users $(1\leq n <P-1)$ departing from the origin experience the trip cost $C^{r}$.
Suppose that there exists a departure time $s_{n+1}^{r}$ such that the $(n+1)$th user departing at this time experiences the same trip cost $C^{r}$ while satisfying $d_{n+1}=d_{n}+m/\mu$.
Then, the following relationship regarding the forecasted cost for $s_{n+1}^{r}$ holds:
\begin{align}
	s_{n+1}^{r}<s_{n+1}
	\quad
	\Rightarrow 
	\quad
	\hat{C}(s_{n+1}^{r}\mid \mathbf{s})\leq C^{r}.
\end{align}
The equality relationship holds only if $V(s_{n+1}^{r})=C^{r}$, which means that $s_{n+1}^{r}=d_{n}+m/\mu$.
\end{lemm}
\begin{prf}
The relationship $s_{n+1}^{r}<s_{n+1}$ implies $s_{n+1}^{r}\in(s_{n},s_{n+1})$.
We thus calculate the forecasted cost referring to the trip costs of $n$th and $(n+1)$th users based on the cases indicated in Eqs.~\eqref{Eq:DefiBetterCongested} and \eqref{Eq:DefiBetterFree}.

Consider first the case where $s_{n+1}\leq d_{n}+m/\mu$, which suggests that $d_{n+1}=d_{n}+m/\mu$ and $s_{n+1}^{r}<d_{n}+m/\mu$.
Substituting this equation into Eq.~\eqref{Eq:TripCost}, we have
\begin{align}
	C_{n+1}=C^{r}-(s_{n+1}-s_{n+1}^{r}).
\end{align}
We then derive the forecasted cost from Eq.~\eqref{Eq:DefiBetterCongested}, as follows:
\begin{align}
	\hat{C}(s_{n+1}^{r}) = 
	C_{n} + \cfrac{C_{n+1} - C_{n}}{s_{n+1} - s_{n} }(s_{n+1}^{r}-s_{n})
	=C^{r}-\cfrac{s_{n+1}-s_{n+1}^{r}}{s_{n+1} - s_{n} }(s_{n+1}^{r}-s_{n})<C^{r}.
\end{align}
This is consistent with the lemma.

Consider next the case where $d_{n}+m/\mu<s_{n+1}$, which means that the forecasted cost is calculated from Eq.~\eqref{Eq:DefiBetterFree}.
First, since the $(n+1)$th user experiences the trip cost $C^{r}$ when arriving at the destination at $d_{n}+m/\mu$ by departing at $s_{n+1}^{r}$, the relationship $V(d_{n}+m/\mu)\leq C^{r}$ holds: otherwise, the schedule cost alone exceeds $C^{r}$.
We also have $V(s_{1})=V(d_{1})=C^{r}$.
By combining these facts and the convexity of the schedule delay cost function, we have
\begin{align}
	V(t)\leq C^{r},\quad \forall t\in(s_{1},d_{n}+m/\mu].\label{Eq:App-ScheduleCostBound}
\end{align}
The equality holds when $V(d_{n}+m/\mu)=C^{r}$ and this holds when $s_{n+1}^{r}=d_{n}+m/\mu$.

We then consider the first case in Eq.~\eqref{Eq:DefiBetterFree}, i.e. the relationship $s_{n+1}^{r} \leq d_{n}$ holds.
In this case, we have $n>1$ since the relationship $d_{1}<s_{2}^{r}$ obviously holds.
This implies that $d_{n}\in(s_{1},d_{n}+m/\mu)$; 
we thus have $V(d_{n})<C^{r}$ from Eq.~\eqref{Eq:App-ScheduleCostBound}.
Substituting this into the equation, we have
\begin{align}
	\hat{C}(s^{r}_{n+1}) =C_{n}+\cfrac{V(d_{n}) - C_{n}}{d_{n} - s_{n}}(s_{n+1}^{r}-s_{n})
	< C^{r}.
\end{align}
This is consistent with the lemma.

We finally consider the second case, i.e. the relationship $d_{n} < s_{n+1}^{r}$ holds.
Since $s_{n+1}^{r}\in (s_{1},d_{n}+m/\mu]$, we have $V(s_{n+1})\leq C^{r}$ from Eq.~\eqref{Eq:App-ScheduleCostBound}.
In addition, the equality relationship only holds only when $V(s_{n+1}^{r}) = V(d_{n}+m/\mu) = C^{r}$.
This is consistent with the lemma. \qed
\end{prf}

\begin{lemm}\label{Lemm:Sn<SnR}
Consider a time profile $\mathbf{s}$ where the first to the $n$th users $(1\leq n <P-1)$ departing from the origin experience the trip cost $C^{r}$ while arriving at the destination at intervals of $m/\mu$.
Suppose that there exists a departure time $s_{n+1}^{r}$ such that the $(n+1)$th user departing at this time experiences the same trip cost $C^{r}$ while satisfying $d_{n+1}=d_{n}+m/\mu$.
Suppose also that $C^{r}\geq \rho$ (i.e. $s_{1}\leq t^{-}$) and $s_{n+1}<s_{n+1}^{r}$.

Then, there exists a better response path from $\mathbf{s}$ to a new profile $\mathbf{s}'$ while satisfying the following three conditions:
\begin{enumerate}
	\item The departure time $s_{n+1}'$ of the $(n+1)$th user in $\mathbf{s}'$ is equal to or later than $s_{n+1}^{r}$.
	\item All better responses are conducted by users from the $(n+1)$th onward in the range of times later than $s_{n}$, 
	i.e. the first to the $n$th users do not conduct better responses, and they are not overtaken by the other users.
	\item When $n<P-1$, there exists a user whose trip cost is equal to or higher than $C^{r}$ among ones departing $(n+1)$th or later in $\mathbf{s}'$.
	\end{enumerate}
\end{lemm}
\begin{prf}
Let $t^{+,r}$ denote the destination arrival time later than the desired arrival time $0$ that satisfies $V(t^{+,r})=C^{r}$.
Regarding this, the relationship $s_{n+1}^{r}\leq t^{+,r}$ holds: otherwise, a user departing at $s_{n+1}^{r}$ arrives at the destination later than $t^{+,r}$ and thus the schedule cost alone exceeds $C^{r}$, which leads to a contradiction.
This also means $s_{n+1}<t^{+,r}$.


First of all, since the $(n+1)$th user departs earlier than $s_{n+1}^{r}$, the destination arrival time is $d_{n}+m/\mu$ in $\mathbf{s}$.
We then have the following relationship:
\begin{align}
C_{n+1}(\mathbf{s})=d_{n+1}-s_{n+1}+V(d_{n+1}) = C^{r} + (s_{n+1}^{r}-s_{n+1})\geq C^{r}+\Delta s\geq \rho+\Delta s.
\end{align}
We then prove conditions 1 and 2 by showing the existence of a departure time $s'>s_{n+1}$ where the forecasted cost is equal to or lower than the current cost $C_{n+1}$.
This is because if such a departure time exists, the $(n+1)$th user can delay the departure time by conducting a better response to $s'$ as long as $s_{n+1}<s_{n+1}^{r}$.
Thus, conditions 1 and 2 are ensured by repeating such better responses.
We also prove condition 3 while proving the existence of such a departure time when $n+1<P$.

\subsubsection*{Case where $n+1=P$}
We first consider the case where $n+1=P$ and show that the $P$th user can conduct a better response to $t^{+,r}$ by calculating the forecasted cost.
When $t^{+,r} > d_{P}$, we have 
\begin{align}
\hat{C}(t^{+,r}) = V(t^{+,r}) = C^{r}.\label{Eq:App-LastTimeForecast}
\end{align}
The user thus can conduct the better response.

When $t^{+,r} \leq d_{P}$, we first have $s_{P}<t^{+,r}$ since $s_{P}<s_{P}^{r}$ and $s_{P}^{r}\leq t^{+,r}$.
This means that the forecasted cost is calculated from Eq.~\eqref{Eq:DefiBetterFree} since the user departing later than $t^{+,r}$ does not exist.
When the second case is applied, we have the same equation as Eq.~\eqref{Eq:App-LastTimeForecast}.
When the first case is applied, we have the following relationship:
\begin{align}
V(d_{P})=V(d_{P-1}+m/\mu)\leq C^{r}.
\end{align}
This is because the arrival time of the $P$th user is $d_{P-1}+m/\mu$ when the user departs at $s_{P}^{r}$ and experiences the trip cost $C^{r}$: the schedule delay cost alone cannot exceed $C^{r}$.
This means that the second term in the first case becomes negative, which implies that $\hat{C}(t^{+,r})<C_{P}$.
The user thus can conduct the better response, which proves the lemma in the case where $n+1=P$.

\subsubsection*{Case where $n+1<P$}
In this case, we further divide the analysis into the following cases: (a) $d_{P}<t^{+,r}$, (b) $t^{+,r}<d_{P}$, (c) $t^{+,r}=d_{P}$ and $C^{r}>\rho$, (d) $t^{+,r}=d_{P}$, $C^{r}=\rho$, $s_{P}\neq d_{P}$ and (e)$t^{+,r}=d_{P}$, $C^{r}=\rho$, $s_{P}= d_{P}$.

We first consider the case where $d_{P}<t^{+,r}$.
This means that $t^{+,r}$ is later than the destination arrival time of the last user.
We thus have $\hat{C}(t^{+,r})=V(t^{+,r})=C^{r}$, which means that the $n+1$st user can conduct a better response to $t^{+,r}$ that is later than $s_{n+1}$.
Moreover, the trip cost of the user is equal to or larger than $C^{r}$ after the better response since the schedule delay cost alone becomes equal to or larger than $C^{r}$.
These prove that the conditions hold in this case.

We next consider the case where $t^{+,r}<d_{P}$.
In this case, the following relationship holds:
\begin{align}
d_{P}-s_{1} > t^{+,r} - t^{-}\geq t^{+}-t^{-}=m(P-1)/\mu.
\end{align}
Combining this and \textbf{Corollary~\ref{Coro:App-FreeFlowRange}}, we see that there exists a departure time $s_{F}$ where $\hat{C}(s_{F})=V(s_{F})\leq \rho$ in $[t^{-},t^{+}]$.
In addition, it is obvious that such a departure time does not exist in $[s_{1},s_{n+1}]$ since they arrive at the destination at the intervals of $m/\mu$ and thus the forecasted cost is calculated from Eq.~\eqref{Eq:DefiBetterCongested}, for congested situations.
This means that $s_{F}$ is later than $s_{n+1}$ and thus the $(n+1)$th user can delay the departure time by conducting a better response to $s_{F}$.
This ensures that the conditions 1 and 2 hold.
In addition, \textbf{Corollary~\ref{Coro:App-FreeFlowRange}} suggests that the departure time of the last user $d_{P}$ does not become earlier by the better response since $d_{P}$ is later than $t^{+,r}$ and thus later than $s_{F}$.
We thus have
\begin{align}
C_{P}\geq V(d_{P})>V(t^{+,r})>C^{r}.
\end{align}
This ensures that the condition 3 holds.

We next consider the case where $t^{+,r}=d_{P}$ and $C^{r}>\rho$.
In this case, we have the following relationship:
\begin{align}
d_{P}-s_{1} > t^{+,r} - t^{-}>t^{+}-t^{-}=m(P-1)/\mu.
\end{align}
Thus, by applying the same logic in the previous case, we can ensure that the conditions 1, 2 and 3 hold.

We further consider the case where $t^{+,r}=d_{P}$, $C^{r}=\rho$ and $s_{P}\neq d_{P}$.
In this case, $d_{P}=t^{+}$ and the forecasted cost is calculated as follows:
\begin{align}
\hat{C}(d_{P})=V(t^{+}) = \rho.
\end{align}
Therefore, the $(n+1)$st can conduct a better response to $t^{+,r}$.
In addition, the trip cost of the user is equal to or larger than $\rho$ after the better response since the schedule delay cost alone becomes equal to or larger than $\rho$.
These prove that all of the conditions hold in this case.

We finally consider the case where $t^{+,r}=d_{P}$, $C^{r}=\rho$ and $s_{P}= d_{P}$, which means that $s_{P}=d_{P}=t^{+}$.
In this case, the $(n+1)$th user cannot conduct a better response to $s_{P}$ since there already exists a user.
We thus show that the user can delay the departure time by conducting a better response to $s_{P}-\Delta s$, which is a departure time later than $s_{P-1}$ and $s_{P-1}^{r}$ $(=s_{P-1}^{e})$.

Since $d_{P}-s_{1}=m(P-1)/\mu$, the interval between the destination arrival times of any consecutive users must be $m/\mu$, as with $\mathbf{s}^{e}$.
Combining this and the fact that $s_{1} = t^{-}$, we see that the destination arrival times of all users are the same as those in $\mathbf{s}^{e}$.
Thus, the trip cost of the $P-1$st user is described as follows:
\begin{align}
    C_{P-1}(\mathbf{s}) = d_{P-1}(\mathbf{s}) - s_{P-1} + V(d_{P-1})  = 
    \rho + (s_{P-1}^{e} - s_{P-1}).
\end{align}
In addition, since the interval between the destination arrival times is $m/\mu$, the forecasted cost for the departure time $s_{P}-\Delta s \in (s_{P-1}, s_{P})$ is calculated from Eq.~\eqref{Eq:DefiBetterCongested} as follows:
\begin{align}
    \hat{C}(s_{P}-\Delta s) 
    &= \rho + 
    \cfrac{ s_{P-1}^{e} - s_{P-1} }{s_{P} - s_{P-1}}\Delta s.
\end{align}
Since $t^{+}>s_{P-1}^{e}$, the relationship $s_{P}>s_{P-1}^{e}$ holds.
This means that the fraction is less than one, and we have
\begin{align}
    \hat{C}(s_{P}-\Delta s)  < \rho + \Delta s.
\end{align}
We thus see that the $(n+1)$th user can delay the departure time by conducting the better response to $s_{P}^{e}-\Delta s$.
In addition, the trip cost of the last user is obviously equal to or larger than $\rho$.

Since the cases are exhausted, the lemma is proved.\qed
	
\end{prf}

\subsection{Proof of Lemma~\ref{Lemm:TheoWAG-1}}
We now prove \textbf{Lemma~\ref{Lemm:TheoWAG-1}} by utilising the above lemmas.
We first consider the case where $s_{n+1}>s_{n+1}^{e}$ holds in the time profile $\mathbf{s}$, in which the first to $n$th users select the equilibrium departure times.
By applying \textbf{Lemma~\ref{Lemm:SnR<Sn}}, we see that $\hat{C}(s_{n+1}^{e})\leq \rho$; 
the equality holds only if $s_{n+1}^{e}=d_{n}+m/\mu$, which means that $n+1=P$.
Then, when $n+1=P$, it is obvious that $C_{P}>\rho$ since the $P$th user departs later than $s_{P}^{e}=t^{+}$.
The user thus can conduct a better response by changing his/her departure time to $s_{P}^{e}$.
When $n+1<P$, \textbf{Lemma~\ref{Lemm:App-FirstLast}} suggests that $C_{P}\geq \rho$.
Since $\hat{C}(s_{n+1}^{e}) < \rho$, the $P$th user can conduct a better response to $s_{n+1}^{e}$.
These ensure the existence of a better response path from $\mathbf{s}$ to another time profile $\mathbf{s}'$ where the first to $(n+1)$st users select the equilibrium departure times.
This proves the lemma in this case.

We next consider the case where $s_{n+1}<s_{n+1}^{e}$ holds in the time profile $\mathbf{s}$
By applying \textbf{Lemma~\ref{Lemm:Sn<SnR}}, we see that there exists a better response path from $\mathbf{s}$ to $\mathbf{s}'$ while satisfying the conditions.
Condition 1 suggests that $s_{n+1}'\geq s_{n+1}^{e}$ in a new time profile $\mathbf{s}'$.
When $s_{n+1}'=s_{n+1}^{e}$, Condition 2 suggests that the first to $(n+1)$th users select the equilibrium departure times in $\mathbf{s}'$; this ensures the existence of the considered better response path.
When $s_{n+1}'>s_{n+1}^{e}$, we can apply \textbf{Lemma~\ref{Lemm:SnR<Sn}} to the new profile $\mathbf{s}'$ and Condition 3 guarantees the existence of a user whose trip cost is equal to or higher than $\rho$.
Then, we can prove the existence of the considered better response path using the same logic as in the case where $s_{n+1}>s_{n+1}^{e}$.
This proves the lemma in this case.

Since the cases are exhausted, the lemma is proved. \qed
\color{black}

\subsection{Proof of Lemma~\ref{Lemm:TheoWAG-2}}
We first consider the case where $s_{1}<t^{-}$; in this case, we show the existence of time profile $\mathbf{s}'$ where the departure time of the first user is later than $\mathbf{s}$.
First, the trip cost of the first user satisfies the following relationship:
\begin{align}
    C_{1}(\mathbf{s}) = -\beta s_{1} > \rho.
\end{align}
It is thus sufficient for us to prove the existence of a departure time where the forecasted cost is equal to or lower than $\rho$.
To this end, we divide the analysis into the following two cases according to the arrival time of the last user $d_{P}(\mathbf{s})$: (a) $d_{P} = s_{1}+m(P-1)/\mu$ and (b) $d_{P}> s_{1}+m(P-1)/\mu$.

In Case (a), the departure time of the last user satisfies
\begin{align}
    d_{P} < t^{-}+\cfrac{m(P-1)}{\mu} = t^{+}.
\end{align}
Thus, the first user can change his/her departure time to $t^{+}$ for which the forecasted cost is $\rho$.
In Case (b), the relationship $d_{P} - s_{1} > m(P-1)/\mu$ holds true.
Thus, from \textbf{Corollary~\ref{Coro:App-FreeFlowRange}}, there exists a departure time where the forecasted cost is equal to or lower than $\rho$.
Therefore, the first user can conduct a better response while delaying the departure time of the first user.

We next consider the case where $t^{-}<s_{1}$; in this case, we show the existence of time profile $\mathbf{s}'$ where the departure time of the first user is equal to $t^{-}$.
Since $t^{-}<s_{1}$, a user can change his/her departure time to $t^{-}$ where the forecasted cost is $\rho$.
Meanwhile, \textbf{Lemma~\ref{Lemm:App-FirstLast}} suggests that the trip cost of the last user is higher than $\rho$.
Therefore, the last user can conduct a better response to $t^{-}$.
These prove the lemma.\qed

\section{Proof of Proposition~\ref{Prop:AsymptoticFirst}}\label{Sec:App-Prop:AsymptoticFirst}
Let $\mathbf{s}$ denote the considered time profile where the first to $n$th users experience the same trip cost $C^{r}$ $(\neq \rho)$ while arriving at the destination at the interval of $m/\mu$ and there does not exist a better response path to a time profile where $(n+1)$st user experiences $C^{r}$.
Such a time profile can be classified into one of the following:
\begin{enumerate}
    \item $n = P$: all users experience the trip cost $C^{r}$.
    
    \item $n < P$ and $V(d_{n}(\mathbf{s})+m/\mu) > C^{r}$: the schedule delay cost of the $n+1$st user exceeds $C^{r}$ when the user arrives at the destination at the interval of $m/\mu$.
    
    \item $n < P$ and $V(d_{n}(\mathbf{s})+m/\mu)) \leq C^{r}$: the schedule delay cost of the $n+1$st user does not exceed $C^{r}$ when the user arrives at the destination at the interval of $m/\mu$.
\end{enumerate}
Below, we prove that the properties of every time profiles are consistent with those in the proposition.
In the analysis, we accord the index of a user with the order of departure as with \textbf{\ref{App:Prf_WAG}}.

\subsection{Case 1}
First of all, this case does not occur when $t^{-} < s_{1}$.
This is because \textbf{Lemma~\ref{Lemm:App-FirstLast}} suggests that $C_{P} > \rho > C^{r}$, which means that the trip cost of the last user is higher than $C^{r}$.
Therefore, it is sufficient to consider the case where $s_{1}< t^{-}$.

Since all users arrive at the destination at the interval of $m/\mu$ in $\mathbf{s}$, the destination arrival time and the trip cost of the last user satisfy the following relationships:
\begin{align}
    d_{P}(\mathbf{s}) = s_{1} + \cfrac{m(P-1)}{\mu} < t^{-} + \cfrac{m(P-1)}{\mu} = t^{+}\quad \Rightarrow\quad V(d_{P}) < \rho <C^{r}.\label{Eq:App-Eventual-1}
\end{align}
This means that the last user does not depart at $d_{P}$ and experiences a positive queueing delay cost.
This is consistent with the property mentioned in the proposition.

\subsection{Case 2}
First of all, this case does not occur when $s_{1}< t^{-}$.
This is because the schedule delay cost of the last user is lower than $C^{r}$ even when all users arrive at the destination at the interval of $m/\mu$, as shown in Eq.~\eqref{Eq:App-Eventual-1}.
Combining this fact and the convexity of the schedule delay cost function, we see that the schedule delay costs of all users do not exceed $C^{r}$ except the first user.
Therefore, the relationship $V(d_{n}(\mathbf{s}) + m/\mu) > C^{r}$ does not hold.

Consider the case where $t^{-} < s_{1}$.
When $V(d_{n}(\mathbf{s})+m/\mu) > C^{r}$ holds, the convexity of the schedule delay cost function implies that 
\begin{align}
    0 < d_{n}(\mathbf{s}) + m/\mu \leq d_{k}(\mathbf{s}),\quad \Rightarrow \quad 
    C_{k}(\mathbf{s})\geq V(d_{k}(\mathbf{s})) > C^{r},\quad 
    \forall k\geq n+1.
\end{align}
This means that the trip cost of the users departing after the $n$th user is higher than $C^{r}$.
This is consistent with the property mentioned in the proposition.

\subsection{Case 3}
In this case, we first establish the following lemma:
\begin{lemm}\label{Lemm:App-ExistSR}
Consider a time profile $\mathbf{s}$ where the first to $n$th user arriving at the destination at the interval of $m/\mu$ $(n<P-1)$ and their trip costs are equal to or higher than $C^{r}$ which is the trip cost of the first user.
Suppose that $V(d_{n}+m/\mu)\leq C^{r}$.
Then, there exists a departure time $s_{n+1}^{r}$ such that the $(n+1)$th user experiences the trip cost $C^{r}$ by departing at that time.
\end{lemm}
\begin{prf}
Suppose that $s_{n+1}^{r}$ satisfies the following relationship: $s_{n}<s_{n+1}^{r}\leq d_{n}+m/\mu$.
When this relationship holds true, the departure time $s_{n+1}^{r}$ is calculated as follows:
\begin{align*}
	&C_{n+1}(s_{n+1}^{r}) = C^{r} = (d_{n+1}-s^{r}_{n+1})+V(d_{n+1})=d_{n}+m/\mu - s^{r}_{n+1} + V(d_{n}+m/\mu)\\
	&\Rightarrow s^{r}_{n+1} = d_{n}+m/\mu + V(d_{n}+m/\mu) - C^{r}
\end{align*}
It is thus sufficient for us to confirm that the $s_{n+1}^{r}$ is included in $(s_{n}, d_{n}+m/\mu]$.

We first have the following relationship:
\begin{align*}
	s^{r}_{n+1} - s_{n}\geq d_{n}+m/\mu + V(d_{n}+m/\mu) - C^{r} - s_{n} 
	=C_{n} - C^{r} + m/\mu + V(d_{n}+m/\mu) - V(d_{n}).
\end{align*}
The supposition suggests $C_{n} - C^{r}\geq 0$.
By applying \textbf{Lemma~\ref{Lemm:App-RelationSchedule}}, we have 
\begin{align*}
	m/\mu + V(d_{n}+m/\mu) - V(d_{n})\geq m(1-\beta)/\mu > 0.
\end{align*}
We thus have $s_{n+1}^{r}-s_{n}>0$.

We next have 
\begin{align*}
	d_{n}+m/\mu - s_{n+1}^{r}
	= C^{r} - V(d_{n}+m/\mu) \geq 0.
\end{align*}
We thus have $d_{n}+m/\mu \geq s_{n+1}^{r}$.
These prove the lemma.
\qed
	
\end{prf}
\noindent This suggests the existence of $s_{n+1}^{r}$ such that the $(n+1)$th user experiences the trip cost $C^{r}$ by departing at that time, in the considered time profile $\mathbf{s}$ in this case.

We then divide the analysis into the following two cases: (a) $s_{1}<t^{-}$ and (b) $t^{-}<s_{1}$.

\subsubsection{Case (a)}
In this case, we prove that there must exist users whose trip costs are different from $C^{r}$, and their trip costs are always lower than $C^{r}$.
We prove this by contradiction:
we suppose that there exist users whose trip costs are equal to or higher than $C^{r}$ among users from the $(n+1)$th onward, and show the existence of a better response path to a new time profile $\mathbf{s}'$ where the $n+1$st user experiences the trip cost $C^{r}$.

Consider first the case where $s_{n+1} > s_{n+1}^{r}$.
By applying \textbf{Lemma~\ref{Lemm:SnR<Sn}}, we derive the forecasted cost for $s_{n+1}^{r}$ as follows:
\begin{align}
	\hat{C}(s_{n+1}^{r})\leq C^{r}.
\end{align}
Then, the equality holds only when $V(s_{n+1}^{r})=C^{r}$.
When the equality holds, it is obvious from the convexity of the schedule delay cost function that $C_{n+1}>C^{r}$; therefore the $(n+1)$th user can conduct a better response to $s_{n+1}^{r}$.
When the equality does not hold, the supposition ensures the existence of a user who can conduct a better response to $s_{n+1}^{r}$.
These ensure the existence of the considered better response path, which contradicts the supposition that such a better response cannot be conducted.

Consider next the case where $s_{n+1} < s_{n+1}^{r}$.
By applying \textbf{Lemma~\ref{Lemm:Sn<SnR}}, we see that there exists a better response path from $\mathbf{s}$ to $\mathbf{s}'$ while satisfying the conditions in the lemma.
Condition 1 suggests that $s_{n+1}'\geq s_{n+1}^{r}$ in a new time profile $\mathbf{s}'$.
When $s_{n+1}'=s_{n+1}^{r}$, Condition 2 suggests that the first to $(n+1)$st users select the equilibrium departure times in $\mathbf{s}'$; this ensures the existence of the considered better response path.
When $s_{n+1}'>s_{n+1}^{e}$, we can apply \textbf{Lemma~\ref{Lemm:SnR<Sn}} to the new profile $\mathbf{s}'$ and prove the existence of the considered better response path using the same logic as in the case where $s_{n+1}>s_{n+1}^{r}$.
These contradict the supposition that such a better response cannot be conducted.

Therefore, there exists a contradiction in any case;
this proves the proposition in this case.

\subsubsection{Case (b)}
In this case, we prove that there must exist users whose trip costs are different from $C^{r}$, and their trip costs are always higher than $C^{r}$.
We prove this by contradiction:
we suppose that there exist users whose trip costs are equal to or lower than $C^{r}$ among users from the $(n+1)$th onward, and show the existence of a better response path to a new time profile $\mathbf{s}'$ where the $n+1$st user experiences the trip cost $C^{r}$.

\color{black}
Consider first the case where $s_{n+1} > s_{n+1}^{r}$.
By applying \textbf{Lemma~\ref{Lemm:SnR<Sn}}, we see that $\hat{C}(s_{n+1}^{r})\leq C^{r} < \rho$.
We also see from \textbf{Lemma~\ref{Lemm:App-FirstLast}} that there exists a user whose trip cost is higher than $\rho$.
Therefore, such a user can conduct a better response to $s_{n+1}^{r}$, which ensures the existence of the considered better response path.

Consider next the case where $s_{n+1} < s_{n+1}^{r}$, which means that $C_{n+1}>C^{r}$.
Suppose that the trip costs of the $(n+1)$th to $(n+k)$th users $(1\leq k < P-n)$ are higher than $C^{r}$ and the cost of the $(n+k+1)$th user is equal to or lower than $C^{r}$.
This suggests that $V(d_{n+k}+m/\mu)\leq C^{r}$.
Combining this fact and \textbf{Lemma~\ref{Lemm:App-ExistSR}}, we confirm the existence of the departure time $s_{n+k}^{r}$ such that the $(n+k)$th user experiences the trip cost $C^{r}$ by departing at that time.
Also, by recursively applying \textbf{Lemma~\ref{Lemm:App-ExistSR}}, we also confirm the existence of such a departure time for all the $(n+1)$th to $(n+k)$th users.

The supposition also suggests that the $(n+1)$th to $(n+k)$th users arrive at the destination at the interval of $m/\mu$, i.e. in congested situations.
This is because if some user experiences a trip cost higher than $C^{r}$ in a free-flow situation, the schedule delay cost alone exceeds $C^{r}$ and users departing later cannot experience a trip cost equal to or lower than $C^{r}$, which contradicts the supposition.
In addition, considering their trip costs, it is obvious that $s_{n+k}<s_{n+k}^{r}<s_{n+k+1}^{r}\leq s_{n+k+1}$, which means that $s_{n+k}^{r}$ exists between $s_{n+k}$ and $s_{n+k+1}$.

Based on these, we prove the following lemma:
\begin{lemm}
	Consider a time profile in this Case (b).
	Suppose that the trip costs of the $(n+1)$th to $(n+k)$th users $(1\leq k < P-n)$ are higher than $C^{r}$ and the cost of the $(n+k+1)$th user is equal to or lower than $C^{r}$.
	Then, the $(n+k)$th user can conduct a better response to $s_{n+k}^{r}$ while maintaining the departure order.
\end{lemm}
\begin{prf}
Consider the case when $d_{n+k+1}-d_{n+k}=m/\mu$.
The forecasted cost for $s_{n+k}^{r}$ is calculated as follows:
\begin{align}
	\hat{C}(s_{n+k}^{r}) = C_{n+k}+ 
	\cfrac{C_{n+k+1} - C_{n+k}}{s_{n+k+1} -  s_{n+k}} ( s_{n+k}^{r} - s_{n+k} )<C_{n+k}.
\end{align}
Therefore, the user can conduct the better response to $s_{n+k}^{r}$.

Consider next the case when $d_{n+k+1}-d_{n+k}>m/\mu$.
Combining the facts that $V(d_{1}) = C^{r}$ and $V(d_{n+k}+m/\mu) \leq C^{r}$ and the convexity of the schedule delay cost function, we have $V(d_{n+k})<C^{r}$.
This means that the forecasted cost calculated by the upper case in Eq.~\eqref{Eq:DefiBetterFree} is lower than $C_{n+k}$.
We also have $V(s_{n+k}^{r})\leq C^{r}$ since $d_{1}< s_{n+k}^{r} < d_{n+k}+m/\mu$ holds true obviously.
Therefore, the user can conduct the better response to $s_{n+k}^{r}$.
These prove the lemma.\qed
\end{prf}
This means that the $(n+k+1)$th user can conduct a better response to $s_{n+k}^{r}$ and the trip cost becomes $C^{r}$ by the better response.
Therefore, by recursively applying the same procedure to users in the backward direction ($k-1, \ldots, n+1$), we see that the $n+1$st user can conduct a better response to $s_{n+1}^{r}$ as the $n+1$st user;
this ensures the existence of the considered better response path.

Therefore, there exists a contradiction in any case.
We thus prove the lemma in this case by contradiction.

Since the cases are exhausted, the proposition is proved.\qed

\color{black}
\section{Proof of Proposition~\ref{Prop:GlobalConvergence}}\label{Sec:App-Prop:GlobalConvergence}

\textbf{Proposition~\ref{Prop:AsymptoticFirst}} and \textbf{Lemma~\ref{Lemm:TheoWAG-2}} ensure that the departure time of the first user $s_{1}$ approaches $t^{-}$ by the adjustment process in Step 4 in the dynamics.
When $s_{1} = t^{-}$, \textbf{Lemma~\ref{Lemm:TheoWAG-1}} ensure that the probability that the time profile transitions to $\mathbf{s}^{e}$ by the fixation process in Steps 2 and 3 is positive.
Since the probability that such a transition occurs becomes one as $\tau\rightarrow \infty$, the better response dynamics almost surely converges to the equilibrium from any initial time profile.\qed
\color{black}

\bibliography{DTC_Atomic}
\bibliographystyle{elsarticle-harv}

\end{document}